\documentclass[a4paper, 11pt, twoside]{article}
\usepackage{amsmath}
\usepackage{amscd}
\usepackage{amssymb}
\usepackage{amsthm}
\usepackage{graphics}
\usepackage{indentfirst}
\usepackage{latexsym}
\usepackage{amsfonts}
\usepackage{multicol}
\usepackage{stmaryrd}
\usepackage{array}
\usepackage{pxfonts}
\usepackage{verbatim}

\newtheoremstyle{mattthm}{}{}{\itshape}{}{\bfseries}{.}{ }{}

\theoremstyle{mattthm}
\newtheorem{lemma}{Lemma}[section]
\newtheorem{propn}[lemma]{Proposition}
\newtheorem{thm}[lemma]{Theorem}
\newtheorem{cory}[lemma]{Corollary}

\newtheoremstyle{mattdef}{}{}{}{}{\bfseries}{.}{ }{}

\theoremstyle{mattdef}
\newtheorem*{rmk}{Remark}
\newtheorem{defn}[lemma]{Definition}
\newtheorem*{eg}{Example}
\newtheorem*{egs}{Examples}
\newtheorem*{rmks}{Remarks}
\newtheorem*{ack}{Acknowledgement}

\begin{document}

\newenvironment{pf}{\noindent\textbf{Proof.}}{\hfill \qedsymbol\newline}
\newenvironment{pfof}[1]{\vspace{\topsep}\noindent\textbf{Proof of {#1}.}}{\hfill \qedsymbol\newline}
\newenvironment{pfenum}{\noindent\textbf{Proof.}\indent\begin{enumerate}\vspace{-\topsep}}{\end{enumerate}\vspace{-\topsep}\hfill \qedsymbol\newline}
\newenvironment{pfnb}{\noindent\textbf{Proof.}}{\newline}

\newcommand\calu{\mathcal{U}}
\newcommand\calf{\mathcal{F}}
\newcommand\bsm{\begin{smallmatrix}}
\newcommand\esm{\end{smallmatrix}}
\newcommand{\rt}[1]{\rotatebox{90}{$#1$}}
\newcommand\la\lambda
\newcommand{\ol}{\overline}
\newcommand{\ul}{\underline}
\newcommand{\lan}{\langle}
\newcommand{\ran}{\rangle}
\newcommand\partn{\mathcal{P}}
\newcommand{\py}[3]{\,_{#1}{#2}_{#3}}
\newcommand{\pyy}[5]{\,_{#1}{#2}_{#3}{#4}_{#5}}
\newcommand{\thmlc}[3]{\textup{\textbf{(\!\! #1 \cite[#3]{#2})}}}
\newcommand{\sss}{\mathfrak{S}_}
\newcommand{\dom}{\trianglerighteqslant}
\newcommand{\doms}{\vartriangleright}
\newcommand{\ndom}{\ntrianglerighteqslant}
\newcommand{\ndoms}{\not\vartriangleright}
\newcommand{\domby}{\trianglelefteqslant}
\newcommand{\domsby}{\vartriangleleft}
\newcommand{\ndomby}{\ntrianglelefteqslant}
\newcommand{\ndomsby}{\not\vartriangleleft}
\newcommand{\subs}[1]{\subsection{#1}}
\newcommand{\nin}{\notin}
\newcommand{\nchar}{\operatorname{char}}
\newcommand{\thmcite}[2]{\textup{\textbf{\cite[#2]{#1}}}\ }
\newcommand\zez{\mathbb{Z}/e\mathbb{Z}}
\newcommand\zepz{\mathbb{Z}/(e+1)\mathbb{Z}}
\newcommand{\bbf}{\mathbb{F}}
\newcommand{\bbc}{\mathbb{C}}
\newcommand{\bbn}{\mathbb{N}}
\newcommand{\bbq}{\mathbb{Q}}
\newcommand{\bbz}{\mathbb{Z}}
\newcommand\zo{\bbn_0}
\newcommand{\gs}{\geqslant}
\newcommand{\ls}{\leqslant}
\newcommand\dw{^\triangle}
\newcommand\wod{^\triangledown}
\newcommand{\hhh}{\mathcal{H}_}
\newcommand{\sect}[1]{\section{#1}}
\newcommand{\ff}{\mathfrak{f}}
\newcommand{\fff}{\mathfrak{F}}
\newcommand\cf{\mathcal{F}}
\newcommand\fkn{\mathfrak{n}}
\newcommand\sx{x}
\newcommand\bra[1]{|#1\ran}
\newcommand\arb[1]{\widehat{\bra{#1}}}
\newcommand\foc[1]{\mathcal{F}_{#1}}
\newcommand{\clam}{\begin{description}\item[\hspace{\leftmargin}Claim.]}
\newcommand{\prof}{\item[\hspace{\leftmargin}Proof.]}
\newcommand{\malc}{\end{description}}
\newcommand\ppmod[1]{\ (\operatorname{mod}\ #1)}
\newcommand\wed\wedge
\newcommand\wede\barwedge
\newcommand\uu[1]{\,\begin{array}{|@{\,}c@{\,}|}\hline #1\\\hline\end{array}\,}
\newcommand{\ux}[1]{\mathfrak{L}_{#1}}
\newcommand\erim{\operatorname{rim}}
\newcommand\mire{\operatorname{rim}'}
\newcommand\mmod{\ \operatorname{Mod}}
\newcommand\cgs\succcurlyeq
\newcommand\cls\preccurlyeq
\newcommand\inc{\mathfrak{A}}%{\operatorname{asc}}

%abacus drawing commands - use a smallmatrix, with \bd for a bead and \nb for a space.
\newcommand\tl{\begin{picture}(8,4)
\put(4,-1){\line(0,1){3}}
\put(4,2){\line(1,0){8}}
\end{picture}}
\newcommand\tll{\begin{picture}(16,4)
\put(4,-1){\line(0,1){3}}
\put(4,2){\line(1,0){8}}
\end{picture}}
\newcommand\tr{\begin{picture}(8,4)
\put(4,-1){\line(0,1){3}}
\put(-4,2){\line(1,0){8}}
\end{picture}}
\newcommand\tm{\begin{picture}(8,4)
\put(4,-1){\line(0,1){3}}
\put(-4,2){\line(1,0){16}}
\end{picture}}
\newcommand\tmm{\begin{picture}(16,4)
\put(4,-1){\line(0,1){3}}
\put(-4,2){\line(1,0){16}}
\end{picture}}
\newcommand{\bd}{\begin{picture}(8,6)
\put(4,-1){\line(0,1){8}}
\put(4,3){\circle*{6}}
\end{picture}}
\newcommand{\nb}{\begin{picture}(8,6)
\put(4,-1){\line(0,1){8}}
\put(3,3){\line(1,0){2}}
\end{picture}}
\newcommand{\vd}{\begin{picture}(8,10)
\put(4,5){\circle*{1}}
\put(4,2){\circle*{1}}
\put(4,8){\circle*{1}}
\end{picture}}
\newcommand{\hd}{\begin{picture}(10,6)
\put(5,3){\circle*{1}}
\put(2,3){\circle*{1}}
\put(8,3){\circle*{1}}
\end{picture}}
\newcommand{\xb}{\begin{picture}(8,6)
\put(4,-4){\line(0,1){14}}
\put(3,3){\line(1,0){2}}
\end{picture}}

%Topmatter
\title{General runner removal and the Mullineux map}
\author{Matthew Fayers\\\normalsize Queen Mary, University of London, Mile End Road, London E1 4NS, U.K.\\\texttt{\normalsize m.fayers@qmul.ac.uk}}
\date{}
\maketitle
\begin{center}
\end{center}
\markboth{Matthew Fayers}{General runner removal and the Mullineux map}
\pagestyle{myheadings}

\begin{abstract}
We prove a new `runner removal theorem' for $q$-decomposition numbers of the level $1$ Fock space of type $A^{(1)}_{e-1}$, generalising earlier theorems of James--Mathas and the author.  By combining this with another theorem relating to the Mullineux map, we show that the problem of finding all $q$-decomposition numbers indexed by partitions of a given weight is a finite computation.
\end{abstract}

%\tableofcontents

\sect{Introduction}

Let $e$ be an integer greater than or equal to $2$, and let $\mathcal{U}$ denote the quantum algebra $U_q(\widehat{\mathfrak{sl}}_e)$ over $\bbq$.  The \emph{level $1$ Fock space} for $\mathcal{U}$ is a $\bbq(q)$-vector space with a standard basis indexed by the set of all partitions.  This has the structure of an integrable $\mathcal{U}$-module, and the submodule generated by the empty partition is isomorphic to the irreducible highest-weight module $L(\Lambda_0)$ for $\mathcal{U}$.  On computing the Lusztig--Kashiwara \emph{canonical basis} for this submodule and expanding with respect to the standard basis, one obtains coefficients $d_{\la\mu}^{e}(q)$ indexed by pairs of partitions $\la$ and $\mu$ with $\mu$ $e$-regular.  These polynomials have become known as `$q$-decomposition numbers' in view of Ariki's proof of the LLT Conjecture, which states that if $\la$ and $\mu$ are partitions of $n$ with $\mu$ $e$-regular, then $d^{e}_{\la\mu}(1)=[S^\la:D^\mu]$, where $S^\la$ and $D^\mu$ denote a Specht module and a simple module for an Iwahori--Hecke algebra at an $e$th root of unity in $\bbc$.  Leclerc and Thibon extended the canonical basis for $L(\Lambda_0)$ to a canonical basis for the whole of the Fock space, yielding $q$-decomposition numbers $d^{e}_{\la\mu}(q)$ for all pairs $(\la,\mu)$ of partitions, and conjectured that when evaluated at $q=1$, these polynomials should give decomposition numbers for appropriate quantised Schur algebras.  This conjecture was proved by Varagnolo and Vasserot \cite{vv}.

$q$-decomposition numbers have been studied extensively in the last ten years or so, with some effort being devoted to finding faster or more enlightening methods for computing the $q$-decomposition numbers.  An important theorem on these lines is the `runner removal theorem' of James and Mathas \cite{jm}, which shows how to equate a $q$-decomposition number $d^{e}_{\la\mu}(q)$ with a `smaller' $q$-decomposition number $d^{e-1}_{\xi\pi}(q)$ under certain conditions based on abacus displays for $\la$ and $\mu$.  This observation that the $q$-decomposition numbers are `independent of $e$' is inherent in Lusztig's famous conjecture for the characters of irreducible modules for reductive algebraic groups, and indeed the James--Mathas theorem admits a very simple proof using the interpretation (due to Goodman and Wenzl \cite{gw} and Varagnolo and Vasserot \cite{vv}) of $q$-decomposition numbers as parabolic Kazhdan--Lusztig polynomials.

In \cite{runrem}, the author proved another runner removal theorem which is in some sense `conjugate' to the James--Mathas theorem, and more recently Chuang and Miyachi \cite{cm} have shown that there are Morita equivalences of $\zeta$-Schur algebras underlying some of these results.  In this paper, we prove a rather stronger runner removal theorem for $q$-decomposition numbers; this includes both the James--Mathas theorem and the author's earlier theorem as special cases.  The way we do this is to define an integer $\ux k(\la)$ associated to a partition $\la$ and an integer $k\in\{0,\dots,e-1\}$, and then to show that if two partitions $\la$ and $\mu$ satisfy $\ux k(\la)=\ux k(\mu)$ for some $k$, then there is a runner which may be removed both abacus displays, resulting in an equality of $q$-decomposition numbers.  The proof of this theorem involves a long calculation using the Leclerc--Thibon algorithm for computing the canonical basis of the Fock space.

In the remainder of the paper, we prove some results which indicate the strength of our main theorem; the main result here is Corollary \ref{finmain}, in which we show that the problem of computing all $q$-decomposition numbers $d^{e}_{\la\mu}(q)$ for partitions $\la$ and $\mu$ of a given $e$-weight \emph{as $e$ varies} is a finite computation.  This requires a further theorem (Theorem \ref{mainmull}), which is the result of a detailed computation describing a relationship between our function $\ux k$ and the Mullineux map.

%xxx
%We now give a brief indication of the layout of this paper.  Section \ref{backsec} contains the basic background material needed in this paper.  In Section \ref{remsec} we describe the procedure of removing a runner from the abacus, and state our main theorems.  In Section \ref{canbas}, we give a detailed account of $q$-decomposition numbers, and prove Theorem \ref{main}.  In Section \ref{mullsec}, we examine the Mullineux map and prove Theorem \ref{mainmull}; this is then used in Section \ref{conseq} to deduce Corollary \ref{finmain}.

\begin{ack}
This research was undertaken while the author was a visiting Postdoctoral Fellow at Massachusetts Institute of Technology, with the support of a Research Fellowship from the Royal Commission for the Exhibition of 1851.  The author is very grateful to M.I.T.\ for its hospitality, and to the 1851 Commission for its generous support.
\end{ack}

\sect{Background}\label{backsec}

%In this section we give some background details on partitions and the abacus, and introduce a slightly unusual version of the dominance order on partitions.

\subsection{Miscellaneous notation}

We begin with some mathematical conventions which might not be considered standard by all readers.
\begin{itemize}
\item
$\zo$ denotes the set of non-negative integers.
\item
If $i$, $j$ and $e$ are integers with $e\gs2$, we write $i\equiv j\ppmod e$ to mean that $i-j$ is divisible by $e$, and we write $i\mmod e$ for the residue of $i$ modulo $e$.
\item
If $I$ and $J$ are multisets of integers, then we write $I\sqcup J$ for the `disjoint union' of $I$ and $J$; that is, the multiset in which the multiplicity of an integer $z$ is the multiplicity of $z$ in $I$ plus the multiplicity of $z$ in $J$.
\item
If $I$ is a multiset of integers and $J$ a set of integers with $J\subseteq I$, then we write $I\setminus J$ to indicate the multiset which consists of $I$ with one copy of each element of $J$ removed.
\item
If $I$ is a multiset of integers and $J$ a set of integers, then we define $|I\cap J|$ to be the number of elements of $I$ \emph{with multiplicity} which are elements of $J$.
\end{itemize}

%.  Details about $q$-decomposition numbers are deferred until Section \ref{canbas}, and similarly a description of the Mullineux map is not given until it is needed in Section \ref{mullsec}.

\subs{Partitions}

As usual, a \emph{partition} of $n$ is a weakly decreasing sequence $\la = (\la_1,\la_2,\dots)$ of non-negative integers whose sum is $n$.  We write $\partn$ for the set of all partitions.  When writing a partition, we usually group equal parts and omit zeroes, so that $(4^2,3,1^3)$ represents the partition $(4,4,3,1,1,1,0,0,\dots)$.  We use $\varnothing$ to denote the \emph{empty partition}, i.e.\ the unique partition of $0$.

If $\la$ is a partition, then the \emph{conjugate partition} $\la'$ is given by
\[\la'_i = \left|\left\{j\gs 1\ \left|\ \la_j\gs i\right.\right\}\right|.\]
If $e$ is an integer greater than $1$ and $\la$ is a partition, then we say that $\la$ is \emph{$e$-regular} if there does not exist $i\gs1$ such that $\la_i=\la_{i+e-1}>0$.  We say that $\la$ is \emph{$e$-restricted} if $\la_i-\la_{i+1}<e$ for all $i$, or equivalently if $\la'$ is $e$-regular.

The \emph{Young diagram} of a partition $\la$ is the set
\[[\la] = \left\{(i,j)\in\bbn^2\ \left|\ j\ls\la_i\right.\right\},\]
whose elements are called \emph{nodes} of $\la$.  A node $(i,j)$ of $\la$ is \emph{removable} if $[\la]\setminus\{(i,j)\}$ is again the Young diagram of a partition, while a pair $(i,j)$ not in $[\la]$ is an \emph{addable node} if $[\la]\cup\{(i,j)\}$ is a Young diagram.  If $(i,j)$ and $(i',j')$ are addable or removable nodes of $\la$, then we say that $(i,j)$ is \emph{above} $(i',j')$ if $i<i'$.

With $e$ fixed as above, we define the \emph{residue} of any node or addable node $(i,j)$ to be $j-i\mmod e$.  Given two partitions $\la$ and $\mu$ and $k\in\{0,\dots,e-1\}$, we write $\la\stackrel k\rightarrow\mu$ to mean that $[\mu]$ is obtained from $[\la]$ by adding an addable node of residue $k$.

\subs{$q$-decomposition numbers}\label{canoni}

Suppose $e$ is an integer greater than or equal to $2$.  The \emph{quantised enveloping algebra} $\calu=U_q(\widehat{\mathfrak{sl}}_e)$ is an associative algebra over $\bbq(q)$ which arises as a deformation of the universal enveloping algebra of the Kac--Moody algebra $\widehat{\mathfrak{sl}}_e$.  $\calu$ has \emph{Chevalley generators} $e_i,f_i,q^h$ for $i\in\zez$ and $h$ lying in the coroot lattice of $\widehat{\mathfrak{sl}}_e$; defining relations between these generators are well-known; for example, see \cite[\S4.1]{llt}.  The \emph{bar involution} is the $\bbq$-linear involution of $\calu$ defined by
\[\ol{e_i}=e_i,\qquad \ol{f_i}=f_i,\qquad \ol{q}=q^{-1},\qquad \ol{q^h}=q^{-h}.\]

In this paper we shall be concerned with a particular $\calu$-module, namely the \emph{level $1$ Fock space} $\calf$.  As a $\bbq(q)$-vector space, this has a `standard' basis $\{\bra\la\mid \la\in\partn\}$ indexed by the set of all partitions.  The $\calu$-module structure on $\calf$ was originally described by Hayashi \cite{hay}, and may be found in \cite[\S4.2]{llt} and many other references; it will suffice for us to describe the action of the generators $f_k$.  If $\la$ is a partition, then $f_k\bra\la$ is a linear combination of vectors $\bra\mu$ indexed by those partitions $\mu$ for which $\la\stackrel k\rightarrow\mu$.  Given such a partition $\mu$, we write $(i,j)$ for the node added to $[\la]$ to obtain $[\mu]$, and define $N(\la,\mu)$ to be the number of addable nodes of $\la$ of residue $k$ above $(i,j)$ minus the number of removable nodes of $\la$ of residue $k$ above $(i,j)$.  Then we have
\[f_k\bra\la = \sum_{\la\stackrel k\rightarrow\mu}q^{N(\la,\mu)}\bra\mu.\]

The Fock space is of particular interest, because the submodule $M$ generated by the vector $\bra\varnothing$ is isomorphic to the irreducible highest-weight representation of $\calu$ with highest weight $\Lambda_0$.  Accordingly, there is a bar involution on $M$, defined by $\ol{\bra\varnothing}=\bra\varnothing$ and $\ol{um} = \ol u\,\ol m$ for $u\in\calu$, $m\in M$.

Leclerc and Thibon \cite{lt} found a way to extend this bar involution to the whole of $\calf$; that is, they defined a bar involution on the whole of $\calf$ which extends the bar involution on $M$ and which is still compatible with the action of $\calu$.  Furthermore, they provided a way to compute the image of a standard basis element $\bra\mu$ under the bar involution, which shows that the image $\ol{\bra\mu}$ of a standard basis element $\bra\mu$ equals $\bra\mu$ plus a linear combination of standard basis elements indexed by partitions dominated by $\mu$ (see below for a definition of the dominance order on partitions).  This unitriangularity property of the bar involution means that one can prove the following.

\begin{thm}\label{bareu}\thmcite{lt}{Theorem 4.1}
For each partition $\mu$, there is a unique vector
\[G(\mu) = \sum_{\la\in\partn}d^e_{\la\mu}(q)\bra\la\in\calf\]
such that:
\begin{itemize}
\item
$\ol{G(\mu)}=G(\mu)$;
\item
$d^e_{\mu\mu}(q)=1$, while $d^e_{\la\mu}(q)$ is a polynomial divisible by $q$ for $\la\neq\mu$.
\end{itemize}
\end{thm}

The vectors $G(\mu)$ form a $\bbq(q)$-basis of $\calf$, which is called the \emph{canonical basis}; this is a \emph{global basis} in the sense of Kashiwara \cite{ka}.

%construct a \emph{canonical basis} (or \emph{lower global basis}) of $\calf$: this is a basis $\{G(\mu)\mid \mu\in\partn\}$ indexed by the set of all partitions, which is uniquely defined by the properties
%\begin{itemize}
%\item
%$\ol{G(\mu)}=G(\mu)$, and
%\item
%when we write $G(\mu)$ as a linear combination $\sum_\la d_{\la\mu}(q)\bra\la$ of standard basis elements, then the coefficient $d_{\mu\mu}(q)$ equals $1$, while for any $\la\neq\mu$ the coefficient $d_{\la\mu}(q)$ is a polynomial divisible by $q$.
%\end{itemize}
This paper is chiefly concerned with computing the transition coefficients $d^e_{\la\mu}(q)$ arising in Theorem \ref{bareu}.  These polynomials are known as \emph{$q$-decomposition numbers}, in view of the following theorem.

\begin{thm}
Suppose $\la$ and $\mu$ are partitions of $n$, and let $\Delta(\la)$ and $L(\mu)$ denote the corresponding Weyl module and simple module for the $\zeta$-Schur algebra $\mathcal{S}_\zeta(n,n)$, where $\zeta$ is a primitive $e$th root of unity in $\bbc$.  Then
\[[\Delta(\la):L(\mu)] = d^e_{\la'\mu'}(1).\]
\end{thm}

This theorem is due to Varagnolo and Vasserot \cite[\S11, Theorem]{vv}; it generalises a version for decomposition numbers of Iwahori--Hecke algebras conjectured by Lascoux, Leclerc and Thibon and proved by Ariki \cite[Theorem 4.4]{ari}.

%The main results of this paper concern the computation of $q$-decomposition numbers.  These are stated in \S\ref{remsec} below; an explicit description of the bar involution is given in \S\ref{canbas}.

\begin{comment}
Suppose $\la$ and $\mu$ are partitions, and $e$ is an integer greater than or equal to $2$.  Then there is a polynomial $d^e_{\la\mu}(q)$, known as a \emph{$q$-decomposition number}, which arises from a canonical basis for a level $1$ Fock space for the quantum algebra $U_q(\widehat{\mathfrak{sl}}_e)$.  (The superscript $e$ does not usually appear in the notation for $q$-decomposition numbers, but is essential to avoid ambiguity in this paper, since we shall be comparing $q$-decomposition numbers for different values of $e$.)  The most important application of $q$-decomposition numbers is the following theorem due to Varagnolo and Vasserot (extending a theorem of Ariki).

\begin{thm}\thmcite{vv}{\S11, Theorem}
Suppose $\la$ and $\mu$ are partitions of $n$, and let $\Delta(\la)$ and $L(\mu)$ denote the corresponding Weyl module and simple module for the $\zeta$-Schur algebra $\mathcal{S}_\zeta(n,n)$, where $\zeta$ is a primitive $e$th root of unity in $\bbc$.  Then
\[[\Delta(\la):L(\mu)] = d^e_{\la'\mu'}(1).\]
\end{thm}

We shall give a more detailed account of $q$-decomposition numbers in Section \ref{canbas}, where we prove our main theorem.
\end{comment}

\subs{The Mullineux map}\label{mullintro}

Fix an integer $e\gs2$.  Given any non-negative integer $n$, there is a bijection $m$ from the set of $e$-regular partitions of $n$ to itself, known as the \emph{Mullineux map}.  This map depends on the value of $e$, and we may write it as $m_e$ if necessary.  This map was introduced by Mullineux in the case where $e$ is a prime, in an attempt to solve the problem of tensoring a simple module for the symmetric group in characteristic $e$ with the one-dimensional signature representation.  Specifically, if $\mu$ is an $e$-regular partition of $n$ and $D^\mu$ is the corresponding simple $\bbf_e\sss n$-module, then the module $D^\mu\otimes\operatorname{sgn}$ is also a simple module, and is therefore labelled by an $e$-regular partition which we denote $M(\mu)$.  Mullineux's conjecture states that $M(\mu)=m(\mu)$ for all $\mu$.  This conjecture was proved by Ford and Kleshchev \cite{fk}, using Kleshchev's alternative combinatorial characterisation \cite{k3} of the map $M$.  Kleshchev's results have since been generalised by Brundan \cite{brun} to Iwahori--Hecke algebras of type $A$ at an $e$th root of unity (where $e$ need no longer be prime), and an analogue of the Mullineux conjecture holds in this context as well.

Our interest in the Mullineux map derives from the following connection with $q$-decomposition numbers; see the next section for the definition of the $e$-weight of a partition.

\begin{propn}\label{lltthm}\thmcite{llt}{Theorem 7.2}
Suppose $\la$ and $\mu$ are partitions with $e$-weight $w$, and that $\mu$ is $e$-regular.  Then
\[d^{e}_{\la'm(\mu)}(q) = q^{w}d^{e}_{\la\mu}(q^{-1}).\]
\end{propn}

Understanding the Mullineux map will be very helpful for us in computing $q$-decomposition numbers.  Our main result concerning the Mullineux map is Theorem \ref{mainmull}; this is proved in Section \ref{mullsec}, where a detailed description of the Mullineux map is given.%, where we prove our second main theorem.  In Section \ref{conseq}, we combine our main theorems with Proposition \ref{lltthm}.ting $q$-decomposition numbers.  Our main result concerning the Mullineux map is Theorem \ref{mainmull}; this is proved in Section \ref{mullsec}, where a detailed description of the Mullineux map is given.%, where we prove our second main theorem.  In Section \ref{conseq}, we combine our main theorems with Proposition \ref{lltthm}.

\subs{The abacus}\label{absec}

Suppose $\la$ is a partition, and $r$ is an integer greater than or equal to $\la'_1$.  For $i=1,\dots,r$ set $\beta_i=\la_i+r-i$.  The integers $\beta_1,\dots,\beta_r$ are distinct, and we refer to the set $\mathfrak{B}_r(\la)=\{\beta_1,\dots,\beta_r\}$ as the \emph{$r$-beta-set} for $\la$.

Now we suppose $e\gs2$, and take an abacus with $e$ vertical runners numbered $0,\dots,e-1$ from left to right.  On runner $i$ we mark positions labelled with the integers $i,i+e,i+2e,\dots$ from the top down.  For example, if $e=4$ then the abacus is marked as follows.
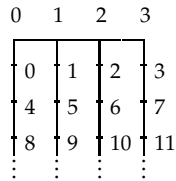
\[\begin{array}{l@{}l@{}l@{}l}
{\ }_0&{\ }_1&{\ }_2&{\ }_3\\[-3pt]\tll&\tmm&\tmm&\tr\\[-3pt]
\xb_{0}&\xb_{1}&\xb_{2}&\xb_{3}\\
\xb_{4}&\xb_{5}&\xb_{6}&\xb_{7}\\
\xb_{8}&\xb_{9}&\xb_{10}&\xb_{11}\\
\vd&\vd&\vd&\vd\\
\end{array}\]
Placing a bead on the abacus at position $\beta_i$ for each $i=1,\dots,r$, we obtain an \emph{abacus display} for $\la$.  In an abacus display, we call a position \emph{occupied} if it contains a bead, and \emph{empty} if it does not, and we we say that position $i$ is \emph{later than} or \emph{after} position $j$ if $i>j$.  For example, we may speak of the `last occupied position' on the abacus, meaning the position $\beta_1$.

Taking an abacus display for $\la$ and sliding all the beads up their runners as far as they will go, we obtain an abacus display for a new partition, which is called the \emph{$e$-core} of $\la$; this partition is independent of the choice of abacus display (i.e.\ the choice of $r$).  The total distance the beads move when we slide them up to obtain the $e$-core of $\la$ is the \emph{$e$-weight} of $\la$.  The $e$-weight and $e$-core are of interest in this paper, because two standard basis vectors $\bra\la$ and $\bra\mu$ lie in the same weight space of $\calf$ if and only if $\la$ and $\mu$ have the same $e$-weight and $e$-core.  (It is unfortunate that the word `weight' is conventionally used in two different ways in this subject; we hope to avoid ambiguity by consistently saying `$e$-weight', reserving `weight' for the Lie-theoretic term.)  Moreover, we have the following statement concerning $q$-decomposition numbers; this essentially says that each canonical basis vector is a weight vector in $\calf$.

\begin{propn}\label{block}\thmcite{lt}{\S4}
Suppose $\la,\mu\in\partn$ with $d^e_{\la\mu}(q)\neq0$.  Then $\la$ and $\mu$ have the same $e$-core and $e$-weight.
\end{propn}

From the definition, we see that two partitions $\la$ and $\mu$ have the same $e$-core if and only if when we take abacus displays for $\la$ and $\mu$ (with the same number of beads on each), there are equal numbers of beads on corresponding runners.  

The abacus formulation provides a way to compare different values of $e$.  Suppose $e$ and $r$ are chosen as above.  Given $k\in\{0,\dots,e-1\}$, we set $d=(r+k)\mmod e$, and we say that a partition $\la$ is \emph{$k$-empty} if all the beads on runner $d$ of the abacus display for $\la$ are as high as possible; that is, there is no $t$ such that position $d+te$ is occupied while $d+(t-1)e$ is empty.  The way we have defined this means that this definition does not depend on the choice of $r$.  The term `$k$-empty' derives from the fact $\la$ is $k$-empty if and only if the $k$th component of the \emph{$e$-quotient} of $\la$ is the empty partition; we refrain from defining and using $e$-quotients in this paper in order to avoid over-complicating notation.

Throughout this paper, if we are given such a triple $(\la,k,r)$, we set $d=(k+r)\mmod e$, and we define $c$ to be the number of beads on runner $d$ of the abacus display for $\la$.  The assertion that $\la$ is $k$-empty means that these beads lie in positions $d,d+e,\dots,d+(c-1)e$.

\begin{eg}
Suppose $e=4$, and $\la$ is the partition $(14,11,9,5,4,1^5)$.  Taking $r=14$, we obtain the following abacus display for $\la$.
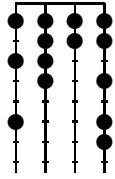
\[\bsm
\tl&\tm&\tm&\tr\\
\bd&\bd&\bd&\bd\\
\nb&\bd&\bd&\bd\\
\bd&\bd&\nb&\nb\\
\nb&\bd&\nb&\bd\\
\nb&\nb&\nb&\nb\\
\bd&\nb&\nb&\bd\\
\nb&\nb&\nb&\bd\\
\nb&\nb&\nb&\nb
\esm\]
We see that $\la$ is $k$-empty for both $k=0$ and $k=3$.  For $k=0$, we have $d=c=2$, while for $k=3$ we have $d=1$, $c=4$.
\end{eg}

\subs{The dominance order}\label{domsec}

Now we introduce a partial order on the set of partitions, which we shall use in place of the usual dominance order.  First we need a partial order on the set of multisets of integers; this is sometimes referred to as the Bruhat order.  Suppose $I = \{i_1,\dots,i_s\}$ and $J=\{j_1,\dots,j_t\}$ are multisets of non-negative integers.  We write $I\cgs J$ if and only if $s=t$ and there is a permutation $\sigma\in\sss s$ such that $i_k\gs j_{\sigma(k)}$ for all $k$.  It is easy to see that $\cgs$ is a partial order.

Now for any finite multiset $B$ of non-negative integers, we define the \emph{$e$-extension} of $B$ to be the multiset $X^e(B)$ of non-negative integers in which the multiplicity of an integer $z$ is
\[\big|B\cap\{z,z+e,z+2e,\dots\}\big|.\]
If $\la$ is a partition and $r$ a large integer, we define the $r$-beta-set $\mathfrak{B}_r(\la)$ for $\la$ as above, and then define the \emph{extended beta-set}
\[\mathfrak{X}^e_r(\la) = X^e(\mathfrak{B}_r(\la)).\]
Given two partitions $\la$ and $\mu$, we say that $\mu$ \emph{dominates} $\la$ (and write $\mu\dom\la$) if $\la$ and $\mu$ have the same $e$-core and $\mathfrak{X}^e_r(\mu)\cgs \mathfrak{X}^e_r(\la)$.  We note that this order does not depend on the choice of $r$.  Indeed, $\mathfrak{X}^e_{r+1}(\la)$ may be obtained from $\mathfrak{X}^e_r(\la)$ by increasing each entry by $1$, and then adding $b$ copies of the integer $0$, where $b$ is the number of beads on runner $0$ of the abacus display for $\la$ with $r+1$ beads.  From this it is easy to see that $\mathfrak{X}^e_r(\mu)\cgs \mathfrak{X}^e_r(\la)$ if and only if $\mathfrak{X}^e_{r+1}(\mu)\cgs \mathfrak{X}^e_{r+1}(\la)$.

We use $\dom$ to denote this order throughout this paper; the usual dominance order will not be used.  Our dominance order depends on the integer $e$, and we may write $\dom_e$ where there is a possibility of ambiguity.

\subs{The Scopes equivalence}

In this section, we briefly recall the Scopes equivalence, as it relates to $q$-decomposition numbers.  Let us define a \emph{block} to be an equivalence class of partitions under the equivalence relation `has the same $e$-core and $e$-weight as'.  In view of Proposition \ref{block}, any non-zero $q$-decomposition number can be regarded as a $q$-decomposition number for (partitions lying in) a particular block.  We define the $e$-core and $e$-weight of a block to be the common $e$-core and $e$-weight of the partitions in that block.  It is easy to see (either combinatorially, or using the fact that $\calf$ has finite-dimensional weight spaces) that any block is finite.

Scopes defined an equivalence relation on the set of blocks of a given $e$-weight, for a given value of $e$.  To describe this, suppose $B$ is a block with $e$-weight $w$ and core $\beta$.  Suppose that for some $k\in\{0,\dots,e-1\}$ $\beta$ has $a$ addable nodes of residue $k$ for some $a\gs w$, and let $\gamma$ be the partition obtained by adding all these addable nodes.  Then $\gamma$ is also an $e$-core; we let $C$ denote the block with $e$-weight $w$ and core $\gamma$.  We say that $B$ and $C$ are \emph{Scopes equivalent}, and we make the Scopes equivalence into an equivalence relation by extending transitively and reflexively.

If $B$ and $C$ are as above, then any $\la\in B$ has exactly $a$ addable nodes of residue $k$ and no removable nodes of this residue.  If we define $\Phi(\la)$ to be the partition obtained by adding these addable nodes, then $\Phi$ is a bijection between $B$ and $C$; these results are proved in \cite[\S2]{js1}.  The condition on addable and removable nodes, together with the action of $f_k$ described in $\S\ref{canoni}$, implies that for any $\la\in B$ we have
\[f_k^{(a)}\bra\la = \bra{\Phi(\la)},\]
where $f^{(a)}_k$ denotes the quantum divided power $f^a_k/[a]!$.  This implies the following.%at the $q$-decomposition numbers are well-behaved under $\Phi$.

\begin{propn}\label{scopq}\thmcite{lm}{Theorem 20}
Let $B$ and $C$ be as above, and take $\la,\mu\in B$.  Then
\[d^e_{\la\mu}(q)=d^e_{\Phi(\la)\Phi(\mu)}(q).\]
\end{propn}

\begin{pf}
Since the $q$-decomposition number $d^e_{\nu\mu}(q)$ is zero unless $\mu$ and $\nu$ have the same $e$-core and $e$-weight, we can write
\[G(\mu) = \sum_{\nu\in B}d^e_{\nu\mu}(q)\bra\nu.\]
Then by the above remarks we have
\[f^{(a)}_kG(\mu) = \sum_{\nu\in B}d^e_{\nu\mu}(q)\bra{\Phi(\nu)}.\]
This vector is invariant under the bar involution (since $G(\mu)$ is, and the bar involution is compatible with the action of $\calu$), and hence by the uniqueness statement in Theorem \ref{bareu} must equal $G(\Phi(\mu))$.
\end{pf}

This gives us an important finiteness result for $q$-decomposition numbers.

\begin{cory}\label{scfin}
Suppose $e,w$ are fixed, and let
\[D^e_w = \left\{\left.d^e_{\la\mu}(q)\ \right|\ \la,\mu\text{ partitions of $e$-weight }w\right\}.\]
Then $D^e_w$ is finite, and there is a finite algorithm to compute it.
\end{cory}

\begin{pf}
For any block $B$, write
\[D_B = \left\{\left.d^e_{\la\mu}(q)\ \right|\ \la,\mu\in B\right\}.\]
Then by Proposition \ref{block}, we have
\[D^e_w = \{0\}\cup\bigcup_BD_B,\]
taking the union over all blocks of $e$-weight $w$.  If $B$ and $C$ are as in Proposition \ref{scopq}, then by that result we have $D_B=D_C$; extending transitively, we get $D_B=D_C$ whenever $B$ and $C$ are Scopes equivalent blocks.  So to compute $D^e_w$ it suffices to consider just one block in each Scopes equivalence class.  Scopes \cite[Theorem 1]{js1} shows that for given $w,e$ there are only finitely many classes (and indicates how to find a representative of each class), and so one only has finitely many blocks to consider.  But any block $B$ is finite, and so the set $D_B$ is finite and may be found with a finite computation.
\end{pf}

\sect{Removing a runner from the abacus -- the main results}\label{remsec}

In this section, we describe the procedure of removing a runner from the abacus, and state our main theorems.

\subs{The runner removal theorems}\label{31}

Suppose $e\gs3$, $k\in\{0,\dots,e-1\}$ and $\la$ is a $k$-empty partition.  Choose a large integer $r$, and let $c,d$ be as defined in \S\ref{absec}.  Construct the abacus display for $\la$ with $r$ beads, and then remove runner $d$.  The resulting configuration will be the abacus display, with $e-1$ runners and $r-c$ beads, for a partition which we denote $\la^{-k}$.  It is a simple exercise to show that the definition of $\la^{-k}$ does not depend on the choice of $r$.

We describe this construction in terms of beta-sets.  Define a function
\[\phi_d:\left\{z\in\zo\mid z\nequiv d\ppmod e\right\}\longrightarrow\zo\]
by setting
\[\phi_d(z) = z-\left\lfloor\frac{z+e-d}e\right\rfloor.\]
Then $\phi_d$ is an order-preserving bijection.  If $\la$ is $k$-empty, then the $r$-beta-set for $\la$ consists of the integers $d,d+e,\dots,d+(c-1)e$, together with some integers $h_1,\dots,h_{r-c}$ not congruent to $d$ modulo $e$.  The set $\{\phi_d(h_1),\dots,\phi_d(h_{r-c})\}$ is then the $(r-c)$-beta-set for $\la^{-k}$.

The idea of removing a runner from the abacus was introduced by James and Mathas, who proved the first `runner removal theorem' for $q$-decomposition numbers; the author subsequently proved a `conjugate' theorem to the James--Mathas theorem.  The idea of these theorems is that if $\la$ and $\mu$ are partitions which are $k$-empty and satisfy some other specified condition, then there is an equality of $q$-decomposition numbers
\[d^{e}_{\la\mu}(q) = d^{e-1}_{\la^{-k}\mu^{-k}}(q).\]
To give precise statements, we suppose that $r,k$ are chosen as above, and set $d=(r+k)\mmod e$.

\begin{thm}\label{runrem1}\thmcite{jm}{Theorem 4.5}
Suppose $e\gs3$ and $\la$ and $\mu$ are $k$-empty partitions with the same $e$-core and $e$-weight.  Suppose that in the $r$-bead abacus displays for each of $\la$ and $\mu$ the last occupied position on runner $d$ is earlier than the first empty position on any runner.  Then
\[d^{e}_{\la\mu}(q) = d^{e-1}_{\la^{-k}\mu^{-k}}(q).\]
\end{thm}

\begin{thm}\label{runrem2}\thmcite{runrem}{Theorem 4.1}
Suppose $e\gs3$ and $\la$ and $\mu$ are $k$-empty partitions with the same $e$-core and $e$-weight.  Suppose that in the $r$-bead abacus displays for each of $\la$ and $\mu$ the first empty position on runner $d$ is later than the last occupied position on any runner.  Then
\[d^{e}_{\la\mu}(q) = d^{e-1}_{\la^{-k}\mu^{-k}}(q).\]
\end{thm}

The theorem we shall prove in this paper generalises both of these theorems; the way we achieve this generality is by finding a condition which is actually a relation between $\la$ and $\mu$ rather than just an absolute condition which both $\la$ and $\mu$ satisfy.  In order to state our theorem, we need to introduce some more notation; the following is the key definition in this paper.

\begin{defn}
Suppose $\la$ is $k$-empty for some $k$; choose a large $r$, and let $c,d$ be as in \S\ref{absec}.  Construct the extended beta-set $\mathfrak{X}^e_r(\la)$ as in \S\ref{domsec}, and then define $\ux k(\la)$ to be the number of elements of $\mathfrak{X}^e_r(\la)$ (with multiplicity) which are greater than $d+ce$.
\end{defn}

It is straightforward to show that the definition of $\ux k(\la)$ does not depend on the choice of $r$.  Now we can state our main theorem.

\vspace{\topsep}

\noindent\hspace{-3pt}\fbox{\parbox{469pt}{\vspace{-\topsep}

\begin{thm}\label{main}
Suppose $e\gs 3$, $\la$ and $\mu$ are partitions with the same $e$-core and $e$-weight, and $k\in\{0,\dots,e-1\}$.  If $\la$ and $\mu$ are $k$-empty and $\ux k(\la)=\ux k(\mu)$, then
\[d^{e}_{\la\mu}(q) = d^{e-1}_{\la^{-k}\mu^{-k}}(q).\]
\end{thm}
\vspace{-\topsep}}}

\begin{rmk}
We comment that the hypothesis `$\la$ and $\mu$ have the same $e$-core and $e$-weight' in Theorems \ref{runrem1}, \ref{runrem2} and \ref{main} is not at all restrictive from the point of view of computing $q$-decomposition numbers, because of Proposition \ref{block}.  We include the hypothesis about the $e$-cores of $\la$ and $\mu$ in order to avoid counterexamples where the abacus displays for $\la$ and $\mu$ have different numbers of beads on runner $d$; and we include the hypothesis concerning the $e$-weights of $\la$ and $\mu$ so that Theorem \ref{main} is actually a generalisation of Theorem \ref{runrem1} (see the remark following Corollary \ref{conux} below).
\end{rmk}

\begin{eg}
Suppose $e=4$, $\la = (7,4,2,1^2)$ and $\mu = (11,2,1^2)$.  Taking $r=9$, we obtain the following abacus displays:
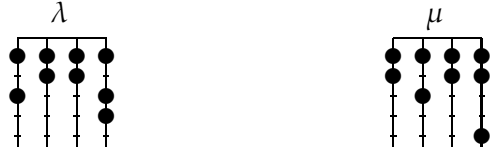
\[\begin{array}c
\la\\
\bsm
\tl&\tm&\tm&\tr\\
\bd&\bd&\bd&\bd\\
\nb&\bd&\bd&\nb\\
\bd&\nb&\nb&\bd\\
\nb&\nb&\nb&\bd\\
\nb&\nb&\nb&\nb
\esm
\end{array}
\qquad\qquad\qquad\qquad
\begin{array}c
\mu\\
\bsm
\tl&\tm&\tm&\tr\\
\bd&\bd&\bd&\bd\\
\bd&\nb&\bd&\bd\\
\nb&\bd&\nb&\nb\\
\nb&\nb&\nb&\nb\\
\nb&\nb&\nb&\bd
\esm
\end{array}\]
We see that $\la$ and $\mu$ both have the same $4$-core and are $1$-empty.  We compute
\begin{align*}
\mathfrak{X}^4_9(\la) &=\{0,0,1,1,2,2,3,3,3,4,5,6,7,7,8,11,11,15\},\\
\mathfrak{X}^4_9(\mu) &=\{0,0,1,1,2,2,3,3,3,4,5,6,7,7,9,11,15,19\}.
\end{align*}
So, taking $c=d=2$, we have $\ux 1(\la)=\ux 1(\mu)=3$, and hence
\[d^4_{\la\mu}(q) = d^3_{\la^{-1}\mu^{-1}}(q),\]
where
\[\la^{-1} = \bsm
\tl&\tm&\tr\\
\bd&\bd&\bd\\
\nb&\bd&\nb\\
\bd&\nb&\bd\\
\nb&\nb&\bd\\
\nb&\nb&\nb
\esm
=(5,3,2,1),\qquad
\mu^{-1} = \bsm
\tl&\tm&\tr\\
\bd&\bd&\bd\\
\bd&\nb&\bd\\
\nb&\bd&\nb\\
\nb&\nb&\nb\\
\nb&\nb&\bd
\esm
=(8,2,1).\]
\end{eg}

Theorem \ref{main} will be proved in Section \ref{canbas}.  We now give our second main theorem, which concerns the relationship between the Mullineux map and the function $\ux k$.  Fix $e\gs2$, and let $m$ denote the Mullineux map (see \S\ref{mullintro}).

\vspace{\topsep}

\noindent\hspace{-3pt}\fbox{\parbox{469pt}{\vspace{-\topsep}

\begin{thm}\label{mainmull}
Suppose $\mu$ is an $e$-regular partition.  If $\mu$ and $m(\mu)'$ are both $k$-empty, then $\ux k(\mu)= \ux k(m(\mu)')$.
\end{thm}
\vspace{-\topsep}}}

\vspace{\topsep}
Theorem \ref{mainmull} will be proved in Section \ref{mullsec}, where a detailed description of the Mullineux map will be given.

\subs{Consequences for the computation of $q$-decomposition numbers}

For the rest of this section, we examine the consequences of Theorems \ref{main} and \ref{mainmull} for the computation of $q$-decomposition numbers.  In order to get to our third main result as quickly as possible, we quote some results which we do not prove until later.

First we need a result which relates $q$-decomposition numbers with the dominance order and the Mullineux map.

\begin{lemma}\label{btwn}
Suppose $\la$ and $\mu$ are partitions with $\mu$ $e$-regular.  If $d^e_{\la\mu}(q)\neq0$, then
\[\mu\dom\la\dom m(\mu)'.\]
\end{lemma}

\begin{pf}
The left-hand inequality is a standard result if $\dom$ is taken to be the usual dominance order; that it holds with our refined dominance order is proved in Proposition \ref{qdecdom} below.  The right-hand inequality follows from this result, together with Proposition \ref{lltthm} and the fact that conjugation of partitions reverses the dominance order (Proposition \ref{condom}).
\end{pf}

Now we prove a useful lemma which seems to be well-known but which the author cannot find in print.

\begin{propn}\label{welk}
Suppose $\la$ and $\mu$ are partitions of $e$-weight $w$.  If $\mu$ is $e$-regular and $\la=m(\mu)'$, then $d^{e}_{\la\mu}(q) = q^w$.  Otherwise, $d^{e}_{\la\mu}(q)$ has degree at most $w-1$.
\end{propn}

\begin{pf}
When $\mu$ is $e$-regular, this is a straightforward consequence of Proposition \ref{lltthm}, since $d^{e}_{\la'm(\mu)}(q)$ is a polynomial which is divisible by $q$ unless $\la'=m(\mu)$.  So we assume that $\mu$ is $e$-singular, and for a contradiction we suppose $d^{e}_{\la\mu}(q)$ has degree at least $w$.

Taking abacus displays for $\la$ and $\mu$, we add a new runner at the left of each display, containing no beads.  Let $\la^+$ and $\mu^+$ be the partitions defined by the resulting displays.  Then by Theorem \ref{runrem1} we have $d^{e}_{\la\mu}(q) = d^{e+1}_{\la^+\mu^+}(q)$.  $\la^+$ and $\mu^+$ obviously have $(e+1)$-weight $w$ and $\mu^+$ is $(e+1)$-regular, and so by the `regular' case of the present proposition we see that $\la^+ = m_{e+1}(\mu^+)'$.  This implies in particular that $\la^+$ is $(e+1)$-restricted, and it is easy to see that this implies that $\la$ is $e$-restricted.  So we can define $\xi = m(\la')$; then $\xi$ is $e$-regular, and $d^e_{\la\xi}(q) = q^w$.  Let $\xi^+$ be the partition obtained by adding an empty runner at the left of the abacus display for $\xi$; then $\xi^+$ is $(e+1)$-regular, and we have $d^{e+1}_{\la^+\xi^+}(q) = d^{e}_{\la\xi}(q) = q^w$.  Applying the `regular' part of the present proposition again, we find that $\la^+ = m_{e+1}(\xi^+)'$.  So
\[\xi^+ = m_{e+1}((\la^+)') = \mu^+,\]
and hence $\xi=\mu$.  But this means that $\mu$ is $e$-regular; contradiction.
\end{pf}

Now we can combine Theorems \ref{main} and \ref{mainmull}.

\begin{propn}\label{kdrop}
Suppose $e\gs3$ and $k\in\{0,\dots,e-1\}$.  Suppose that $\mu$ is an $e$-regular partition such that both $\mu$ and $m(\mu)'$ are $k$-empty.  Then
\begin{enumerate}
\item
$\mu^{-k}$ is $(e-1)$-regular, with $m_{e-1}(\mu^{-k})'=(m(\mu)')^{-k}$, and
\item
for every partition $\la$ with $d^{e}_{\la\mu}(q)\neq0$, we have $d^{e}_{\la\mu}(q)=d^{e-1}_{\la^{-k}\mu^{-k}}(q)$.
\end{enumerate}
\end{propn}

\begin{pf}
We prove (2) first.  Since $d^{e}_{\la\mu}(q)\neq0$, we have $\mu\dom\la\dom m(\mu)'$ by Lemma \ref{btwn}, so by Lemma \ref{sandwich} below and Theorem \ref{mainmull} $\la$ is $k$-empty and $\ux k(\la)=\ux k(\mu)$.  Hence by Theorem \ref{main} we have $d^{e-1}_{\la^{-k}\mu^{-k}}(q)=d^{e}_{\la\mu}(q)$.

Now we prove (1).  Putting $\la=m(\mu)'$ in (2), we have $d^{e-1}_{\la^{-k}\mu^{-k}}(q)=q^w$, where $w$ is the $e$-weight of $\mu$ (and therefore the $(e-1)$-weight of $\mu^{-k}$).  Now (1) follows from Proposition \ref{welk}.
\end{pf}

Now we can show that every $q$-decomposition number which occurs for partitions of $e$-weight $w$ occurs with $e\ls 2w$.

\begin{cory}\label{finite}
Suppose $e>2w>0$.  If $\la,\mu$ are partitions with $e$-weight $w$ and with $d^e_{\la\mu}(q)\neq0$, then there are partitions $\xi,\rho$ of $(2w)$-weight $w$ such that $d^{e}_{\la\mu}(q) = d^{2w}_{\xi\rho}(q)$.
\end{cory}

\begin{pf}
We begin with the case where $\mu$ is $e$-regular.  There are at least $e-w$ values of $k$ for which $\mu$ is $k$-empty, and at least $e-w$ values for which $m(\mu)'$ is $k$-empty.  Since $e>2w$, there must therefore be some $k$ such that $\mu$ and $m(\mu)'$ are both $k$-empty.  By Proposition \ref{kdrop}, $\mu^{-k}$ is $(e-1)$-regular and $d^{e}_{\la\mu}(q)$ equals $d^{e-1}_{\la^{-k}\mu^{-k}}(q)$.  By induction on $e$, this equals $d^{2w}_{\xi\rho}(q)$ for some $\xi,\rho$ of $(2w)$-weight $w$.

Now we consider the case where $\mu$ is $e$-singular.  In this case, we take abacus displays for $\la$ and $\mu$, and add an empty runner at the left of each display.  The resulting partitions $\la^+$ and $\mu^+$ both have $(e+1)$-weight $w$, and $\mu^+$ is $(e+1)$-regular, and we have $d^{e}_{\la\mu}(q)=d^{e+1}_{\la^+\mu^+}(q)$.  So we may apply the `regular' case of the present proposition, replacing $\la,\mu,e$ with $\la^+,\mu^+,e+1$.
\end{pf}

This yields our third main result.

\vspace{\topsep}
\noindent\hspace{-3pt}\fbox{\parbox{469pt}{\vspace{-\topsep}

\begin{cory}\label{finmain}
Fix $w>0$.  Then the set
\[\left\{\left.d^{e}_{\la\mu}(q)\ \right|\ e\gs2,\ \la,\mu\text{ partitions of $e$-weight $w$}\right\}\]
is finite, and there is a finite algorithm to compute it.
\end{cory}
\vspace{-\topsep}}}

\vspace{\topsep}

\begin{pf}
By Corollary \ref{finite}, the given set equals
\[\left\{\left.d^{e}_{\la\mu}(q)\ \right|\ 2\ls e\ls 2w,\ \la,\mu\text{ partitions of $e$-weight $w$}\right\}=D^2_w\cup D^3_w\cup\dots\cup D^{2w}_w,\]
where $D^e_w$ is given in Corollary \ref{scfin}.  Now Corollary \ref{scfin} gives the result.
\end{pf}

\begin{comment}
So it suffices to show that for any $e$ the set
\[D_e=\left\{\left.d^{e}_{\la\mu}(q)\ \right|\ \la,\mu\text{ partitions of $e$-weight $w$}\right\}\]
is finite, and may be computed with a finite algorithm.  But this is well known to be the case: if we fix $e$ and define a \emph{block} to be an equivalence class of partitions under the relation `has the same $e$-weight and $e$-core as', then there is an equivalence relation (the \emph{Scopes equivalence}) on the set of blocks, with the following properties:
\begin{itemize}
\item
if $B$ and $C$ are Scopes equivalent blocks, then there is a bijection $\Phi:B\rightarrow C$ such that for any $\la,\mu\in B$ we have $d^{e}_{\la\mu}(q)=d^{e}_{\Phi(\la)\Phi(\mu)}(q)$;
\item
there are only finitely many Scopes equivalence classes of blocks of a given $e$-weight.
\end{itemize}
The Scopes equivalence was introduced for representation theory of symmetric groups by Scopes \cite{js1}, and extended to Iwahori--Hecke algebras by Jost \cite{jost}.  As a consequence of these results, we see that to compute $D_e$ we need only compute the $q$-decomposition numbers for finitely many blocks, and this is a finite computation.
\end{comment}

As an application of Corollary \ref{finmain}, consider the case $w=3$.  In this case, one can check that each $q$-decomposition number $d^{e}_{\la\mu}(q)$ is always zero or a monic monomial; specifically, it always equals $0$, $1$, $q$, $q^2$ or $q^3$.  This is a significant part of a more general theorem (namely, that all decomposition numbers for weight three blocks of Iwahori--Hecke algebras of type $A$ in characteristic at least $5$ are either $0$ or $1$) which presented great difficulties for several years, until finally proved by the author \cite{mf3}.  Using the above results, more than half of the proof in \cite{mf3} may be replaced by a short computer calculation.

\sect{Combinatorial results}\label{someutil}

In this section, we examine the combinatorics of dominance and runner removal, and prove some simple results which will be useful in the rest of the paper.

\subs{Conjugation and the dominance order}\label{conj}

First we examine how conjugation of partitions relates to the abacus and the dominance order.  Suppose $\la$ is a partition and $r,s$ are large integers.  Let $\mathfrak{B}_r(\la) = \{\beta_1,\dots,\beta_r\}$ be the $r$-beta-set for $\la$, and let $\mathfrak{B}_s(\la')=\{\gamma_1,\dots,\gamma_s\}$ be the $s$-beta-set for $\la'_1$.  The following relationship between these beta-sets is well-known and easy to prove.

\begin{lemma}\label{conbeta}
The set $\{0,\dots,r+s-1\}$ is the disjoint union of $\{\beta_1,\dots,\beta_r\}$ and $\{r+s-1-\gamma_1,\dots,r+s-1-\gamma_s\}$.
\end{lemma}

To express this result in terms of abacus displays, suppose that $r+s\equiv 0\ppmod e$.  Then the $s$-bead abacus display for $\la'$ may be obtained from the $r$-bead abacus display for $\la$ by truncating the diagram after position $r+s-1$, rotating through $180^\circ$, replacing each bead with an empty space, and replacing each empty space with a bead.  As a consequence, we see that the $e$-core of $\la'$ is the conjugate of the $e$-core of $\la$; hence two partitions $\la$ and $\mu$ have the same $e$-core if and only if $\la'$ and $\mu'$ have the same $e$-core.

\begin{eg}
Suppose $e=3$, and $\la=(9,7,4^2,2,1^3)$, so that $\la'=(8,5,4^2,2^3,1^2)$.  Abacus displays for $\la$ and $\la'$, with $11$ and $13$ beads respectively, are as follows, and Lemma \ref{conbeta} can easily be checked.

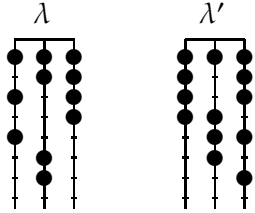
\[\begin{array}c
\la\\
\bsm
\tl&\tm&\tr\\
\bd&\bd&\bd\\
\nb&\bd&\bd\\
\bd&\nb&\bd\\
\nb&\nb&\bd\\
\bd&\nb&\nb\\
\nb&\bd&\nb\\
\nb&\bd&\nb\\
\nb&\nb&\nb
\esm
\end{array}
\qquad
\begin{array}c
\la'\\
\bsm
\tl&\tm&\tr\\
\bd&\bd&\bd\\
\bd&\nb&\bd\\
\bd&\nb&\bd\\
\bd&\bd&\nb\\
\nb&\bd&\bd\\
\nb&\bd&\nb\\
\nb&\nb&\bd\\
\nb&\nb&\nb
\esm
\end{array}
\]
\end{eg}

One consequence of Lemma \ref{conbeta} is the fact that conjugation reverses the dominance order.

\begin{propn}\label{condom}
Suppose $\la$ and $\mu$ are partitions.  Then $\mu\dom\la$ if and only if $\la'\dom\mu'$.
\end{propn}

To prove this, we give an alternative characterisation of the dominance order.  Given a partition $\la$ and a large integer $r$, we define the \emph{weight multiset} $\mathfrak{W}^e_r(\la)$ as follows.  We construct the abacus display for $\la$ with $r$ beads, and slide the beads up their runners to obtain an abacus display for the $e$-core of $\la$.  Each time we slide a bead up one space, from position $b$ to position $b-e$, say, we add a copy of the integer $b$ to the multiset $\mathfrak{W}^e_r(\la)$.

Note that when we slide a bead from position $b$ to position $b-e$, we remove a copy of $b$ from the extended beta-set $\mathfrak{X}^e_r(\la)$.  So if we let $\kappa$ denote the $e$-core of $\la$, then we see that
\[\mathfrak{X}^e_r(\la) = \mathfrak{X}^e_r(\kappa)\sqcup \mathfrak{W}^e_r(\la).\]

\begin{eg}
Suppose $e=2$ and $\la=(3^3)$.  Then $\kappa=(1)$, and taking $r=5$ we get the following abacus displays for $\la$ and $\kappa$.
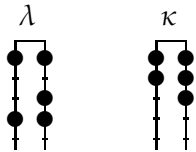
\[
\begin{array}c
\la\\
\bsm
\tl&\tr\\
\bd&\bd\\
\nb&\nb\\
\nb&\bd\\
\bd&\bd\\
\nb&\nb
\esm
\end{array}
\qquad
\begin{array}c
\kappa\\
\bsm
\tl&\tr\\
\bd&\bd\\
\bd&\bd\\
\nb&\bd\\
\nb&\nb\\
\nb&\nb
\esm
\end{array}
\]
We can compute
\begin{align*}
\mathfrak{X}^2_5(\kappa) &= \{0,0,1,1,1,2,3,3,5\}.\\
\mathfrak{W}^2_5(\la) &= \{4,5,6,7\},\\
\mathfrak{X}^2_5(\la) &= \{0,0,1,1,1,2,3,3,4,5,5,6,7\}.
\end{align*}
\end{eg}

Now it is easy to see the following.

\begin{lemma}\label{altdom}
Suppose $\la$ and $\mu$ are partitions with the same $e$-core, and $r$ is a large integer.  Then $\mu\dom\la$ if and only if $\mathfrak{W}^e_r(\mu)\cgs \mathfrak{W}^e_r(\la)$.
\end{lemma}

Notice that the cardinality of $\mathfrak{W}^e_r(\la)$ is the $e$-weight of $\la$; so Lemma \ref{altdom} implies in particular that if $\mu\dom\la$ then $\la$ and $\mu$ have the same $e$-weight.

\begin{pfof}{Proposition \ref{condom}}
Assume $\la$ and $\mu$ have the same $e$-core, and choose large integers $r,s$.  Define $\hat a=r+s+e-1-a$ for any integer $a$.  By Lemma \ref{altdom}, it suffices to show that if $\mathfrak{W}^e_r(\mu)\cgs \mathfrak{W}^e_r(\la)$ then $\mathfrak{W}^e_s(\la')\cgs \mathfrak{W}^e_s(\mu')$.  But by Lemma \ref{conbeta}, we see that moving a bead from position $b$ to position $b-e$ in the $r$-bead abacus display for $\la$ corresponds to moving a bead from position $\hat b$ to position $\hat b-e$ in the $s$-bead abacus display for $\la'$.  So if $\mathfrak{W}^e_r(\la)=\{b_1,\dots,b_w\}$ and $\mathfrak{W}^e_r(\mu)=\{c_1,\dots,c_w\}$, then we have
\[\mathfrak{W}^e_s(\la') = \{\hat b_1,\dots,\hat b_w\},\qquad \mathfrak{W}^e_s(\mu') = \{\hat c_1,\dots,\hat c_w\}.\]
Now if $\sigma\in\sss w$ is such that $c_i\gs b_{\sigma(i)}$ for all $i$, then $\hat b_i\gs \hat c_{\sigma^{-1}(i)}$ for all $i$.
\end{pfof}

\subs{Runner removal and the dominance order}

\begin{lemma}\label{rembij}
Suppose $\la$ is a $k$-empty partition.  Then the map $\xi\mapsto\xi^{-k}$ is a bijection from the set
\[\big\{\xi\in\partn\ \big|\ \xi\text{ is $k$-empty and has the same $e$-core as $\la$}\big\}\]
to the set
\[\big\{\pi\in\partn\ \big|\ \pi\text{ has the same $(e-1)$-core as $\la^{-k}$}\big\}.\]
Furthermore, if $\xi$ is a $k$-empty partition with the same $e$-core as $\la$, then we have $\xi\dom_e\la$ if and only if $\xi^{-k}\dom_{e-1}\la^{-k}$.
\end{lemma}

\begin{pf}
The first statement is obvious from the construction.  For the second statement, we use the characterisation of the dominance order in terms of weight multisets from \S\ref{conj}, and observe that if $\mathfrak{W}^e_r(\la)=\{b_1,\dots,b_w\}$, then $b_1,\dots,b_w\nequiv d\ppmod e$ and
\[\mathfrak{W}^{e-1}_{r-c}(\la^{-k}) = \left\{\phi_d(b_1),\dots,\phi_d(b_w)\right\}.\]
A similar statement applies to $\mathfrak{W}^{e-1}_{r-c}(\xi)$, and the fact that $\phi_d$ is order-preserving gives the result.
\end{pf}

\begin{lemma}\label{uxdom}
Suppose $\la$ and $\mu$ are $k$-empty partitions satisfying $\mu\dom\la$.  Then $\ux k(\mu)\gs \ux k(\la)$.
\end{lemma}

\begin{pf}
This is immediate from the definitions.
\end{pf}

\begin{eg}
We now demonstrate why it is important that we use our coarse version of the dominance order in this paper.  Suppose $e=9$, $\mu=(9,5^2,2^4)$ and $\la=(7^2,3^3,1^4)$.  Taking $r=9$, we get beta-sets
\begin{align*}
\mathfrak{B}_9(\mu) &= \{17,12,11,7,6,5,4,1,0\},\\
\mathfrak{B}_9(\la) &= \{15,14,9,8,7,4,3,2,1\},
\end{align*}
giving abacus displays as follows.
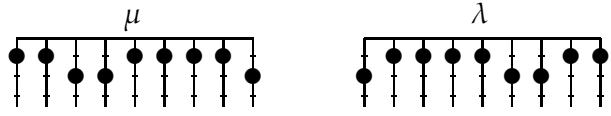
\[\begin{array}c\mu\\\begin{smallmatrix}
\tl&\tm&\tm&\tm&\tm&\tm&\tm&\tm&\tr\\
\bd&\bd&\nb&\nb&\bd&\bd&\bd&\bd&\nb\\
\nb&\nb&\bd&\bd&\nb&\nb&\nb&\nb&\bd\\
\nb&\nb&\nb&\nb&\nb&\nb&\nb&\nb&\nb
\end{smallmatrix}\end{array}
\qquad
\begin{array}c\la\\\begin{smallmatrix}
\tl&\tm&\tm&\tm&\tm&\tm&\tm&\tm&\tr\\
\nb&\bd&\bd&\bd&\bd&\nb&\nb&\bd&\bd\\
\bd&\nb&\nb&\nb&\nb&\bd&\bd&\nb&\nb\\
\nb&\nb&\nb&\nb&\nb&\nb&\nb&\nb&\nb
\end{smallmatrix}\end{array}\]
We see that $\la$ and $\mu$ are $4$-empty.  The extended beta-sets are
\begin{align*}
\mathfrak{X}^9_9(\mu) &= \{0,1,2,3,4,5,6,7,8,11,12,17\},\\ \mathfrak{X}^9_9(\la)&=\{0,1,2,3,4,5,6,7,8,9,14,15\},
\end{align*}
from which we see that $\ux 4(\mu)=1$, $\ux 4(\la)=2$ and $\mu\ndom\la$.  But it is easy to check that $\mu\dom\la$ in the usual dominance order.  So in Lemma \ref{uxdom} (and in everything that follows from it) we need to use our dominance order.
\end{eg}

\begin{lemma}\label{sandwich}
Suppose $\la$ and $\mu$ are $k$-empty partitions satisfying $\mu\dom\la$ and $\ux k(\la)=\ux k(\mu)$.  Then any partition $\xi$ such that $\mu\dom\xi\dom\la$ is $k$-empty, with $\ux k(\xi)=\ux k(\mu)$.  Furthermore, the map $\xi\mapsto\xi^{-k}$ defines a bijection between the sets
\[\{\xi\in\partn\mid \mu\dom_e\xi\dom_e\la\}\]
and
\[\{\pi\in\partn\mid \mu^{-k}\dom_{e-1}\pi\dom_{e-1}\la^{-k}\}.\]
\end{lemma}

\begin{pf}
It suffices to prove that $\xi$ is $k$-empty; the other statements then follow from Lemmata \ref{rembij} and \ref{uxdom}.  So suppose for a contradiction that $\xi$ is not $k$-empty.  The fact that $\mu\dom\xi\dom\la$ implies that $\la,\xi,\mu$ all have the same $e$-core.  So if we take a large integer $r$ and define $c,d$ as in \S\ref{absec}, then $c$ is the number of beads on runner $d$ in the abacus display for each of $\la,\xi,\mu$.  The fact that $\xi$ is not $k$-empty implies that in the abacus display for $\xi$, there is a bead at position $d+te$ for some $t\gs c$.  In particular, this means that the extended beta-set $\mathfrak{X}^e_r(\xi)$ contains the integer $d+ce$.

Write
\begin{align*}
\mathfrak{X}^e_r(\la) &= \{l_1,\dots,l_s\},\\
\mathfrak{X}^e_r(\xi) &= \{x_1,\dots,x_s\},\\
\mathfrak{X}^e_r(\mu) &= \{m_1,\dots,m_s\},
\end{align*}
choosing the ordering so that $l_i\ls x_i\ls m_i$ for each $i$, and $x_1=d+ce$.  Recall that
\[\ux k(\la) = \big|\{i\mid l_i>d+ce\}\big|\]
and similarly for $\mu$.  Now given our choice of numbering, we see that
\[\{i\mid l_i>d+ce\}\subset\{i\mid m_i>d+ce\};\]
the inclusion is strict because $l_1<d+ce<m_1$.  So $\ux k(\la)<\ux k(\mu)$, which is a contradiction.
\end{pf}

\subs{Runner removal and conjugation}
%Next, we compare runner removal and the function $\ux k$ for two conjugate partitions.

%Now we give an equivalent definition of the function $\ux k$, which will be useful in computations.  The proof is an easy exercise.

Now we examine the relationship between runner removal and conjugation.  We begin with a lemma which gives an alternative way to compute $\ux k(\la)$; the proof of this is an easy exercise.

\begin{lemma}\label{altux}
Suppose $\la$ is a $k$-empty partition and $r$ is a large integer.  Let $c,d$ be as defined in \S\ref{absec}, and for any integer $\beta$ let
\[\mathring{\beta} = \begin{cases}
\ \left\lfloor\dfrac{\beta-d-(c-1)e}e\right\rfloor & (\beta\gs d+(c-1)e)\\[9pt]
\qquad\ \ \,\quad0 & (\beta< d+(c-1)e).
\end{cases}\]
Then
\[\ux k(\la) = \sum_{\beta\in\mathfrak{B}_r(\la)}\mathring\beta.\]
\end{lemma}

\begin{lemma}\label{remcon}
Suppose $\la$ is a partition with $e$-core $\kappa$ and $e$-weight $w$, and $k\in\{0,\dots,e-1\}$.  Then $\la$ is $k$-empty if and only if $\la'$ is $(e-1-k)$-empty.  If this is the case, then
\[\ux k(\la)+\ux{e-1-k}(\la') = \ux k(\kappa)+\ux{e-1-k}(\kappa')+w.\]
\end{lemma}

\begin{pf}
We construct abacus displays for $\la$ and $\la'$ using $r$ beads and $s$ beads respectively, where for convenience we choose $r$ and $s$ such that $e\mid r+s$.  Let $c,d$ be as defined in \S\ref{absec}, and set
\[\check d = e-1-d,\qquad \check c = \frac{r+s}e-c.\]
Then $\check d = (s+(e-1-k))\mmod e$, and by Lemma \ref{conbeta} there are $\check c$ beads on runner $\check d$ of the abacus display for $\la'$.  $\la'$ is $(e-1-k)$-empty if and only if these beads are in positions $\check d,\check d+e,\dots,\check d+(\check c-1)e$.  By Lemma \ref{conbeta}, this is equivalent to the condition that the beads on runner $d$ of the abacus display for $\la$ are in positions $d,d+e,\dots,d+(c-1)e$, i.e.\ $\la$ is $k$-empty.

Now suppose $\la$ is $k$-empty.  Assuming $w>0$, we may slide a bead one space up its runner, to obtain a new partition $\xi$.  $\xi$ has $e$-core $\kappa$ and $e$-weight $w-1$ and is $k$-empty, so by induction on $w$ it suffices to show that $\ux k(\la)+\ux{e-1-k}(\la') = \ux k(\xi)+\ux{e-1-k}(\xi')+1$.

Suppose that to obtain $\xi$ from $\la$ we move a bead from position $b$ to position $b-e$.  Then we have $\mathfrak{X}^e_r(\xi)=\mathfrak{X}^e_r(\la)\setminus\{b\}$, and therefore
\[\ux k(\xi) = \begin{cases}
\ux k(\la) & (b<d+ce)\\
\ux k(\la)-1 & (b>d+ce).
\end{cases}\]
By Lemma \ref{conbeta}, an abacus display for $\xi'$ is obtained from an abacus display for $\la'$ by moving a bead from position $r+s-1-b+e$ to position $r+s-1-b$, so
\[\mathfrak{X}^e_s(\xi') = \mathfrak{X}^e_s(\la')\setminus\{r+s-1-b+e\}.\]
By definition $\ux{e-1-k}(\la')$ is the number of elements of $\mathfrak{X}^e_s(\la')$ greater than $\check d+\check ce$, so
\[\ux{e-1-k}(\xi') = \begin{cases}
\ux{e-1-k}(\la') & (r+s-1-b+e<\check d+\check ce)\\
\ux{e-1-k}(\la')-1 & (r+s-1-b+e>\check d+\check ce).
\end{cases}\]
Retracing the definitions gives $r+s-1-b+e<\check d+\check ce$ if and only if $b>d+ce$, and the result follows.
\end{pf}

Now the following is immediate.

\begin{cory}\label{conux}
Suppose $\la$ and $\mu$ are $k$-empty partitions with the same $e$-core and the same $e$-weight.  Then
\[\ux k(\la)+\ux{e-1-k}(\la') = \ux k(\mu)+\ux{e-1-k}(\mu').\]
In particular, $\ux k(\la)=\ux k(\mu)$ if and only if $\ux{e-1-k}(\la')=\ux{e-1-k}(\mu')$.
\end{cory}

\begin{rmk}
We can now show that Theorem \ref{main} is a generalisation of Theorems \ref{runrem1} and \ref{runrem2}.  To see this, it suffices to show that if $\la$ and $\mu$ satisfy the conditions of one of these theorems then $\ux k(\la)=\ux k(\mu)$.  For Theorem \ref{runrem2} this is easy, since by Lemma \ref{altux} the condition in that theorem implies that $\ux k(\la)=0=\ux k(\mu)$.  For Theorem \ref{runrem1}, we note that if $\la$ and $\mu$ satisfy the given conditions, then (by Lemma \ref{conbeta}) the partitions $\la'$ and $\mu'$ satisfy the hypotheses of Theorem \ref{runrem2} (with $k$ replaced by $e-1-k$).  Hence $\ux{e-1-k}(\la')=0=\ux{e-1-k}(\mu')$, and so by Corollary \ref{conux} we have $\ux k(\la)=\ux k(\mu)$.
\end{rmk}

\subs{An alternative characterisation}\label{analter}

The next result in this section gives an alternative characterisation of the relation $\ux k(\la)=\ux k(\mu)$, for two partitions $\la$ and $\mu$.  Given a large integer $r$, define $c,d$ as in \S\ref{absec}, and set
\[n_{r,k}(\la) = \Big|\big\{a,b\in \mathfrak{B}_r(\la)\ \big|\ a<b,\ a\equiv d\nequiv b\ppmod e\big\}\Big|.\]

\begin{lemma}\label{nrk}
Suppose $\la$ is a $k$-empty partition, with $e$-weight $w$ and $e$-core $\kappa$.  If $r$ is a large integer, then
\[\ux k(\la)+n_{r,k}(\la) = \ux k(\kappa)+n_{r,k}(\kappa)+w.\]
\end{lemma}

\begin{pf}
If $\la=\kappa$ then there is nothing to prove, so we assume otherwise.  Then in the abacus display for $\la$ with $r$ beads, we may slide a bead up its runner, from position $b$ to position $b-e$, say, to obtain a new partition $\xi$ with $e$-weight $w-1$.  By induction on $w$, it suffices to prove that
\[\ux k(\la)+n_{r,k}(\la) = \ux k(\xi)+n_{r,k}(\xi)+1.\]
The $r$-beta-set $\mathfrak{B}_r(\la)$ consists of the integers $d,d+e,\dots,d+(c-1)e$, together with $r-c$ integers not congruent to $d$ modulo $e$; the same statement applies to $\mathfrak{B}_r(\xi)$.  Write
\begin{align*}
N_{r,k}(\la) &= \big\{a,b\in\mathfrak{B}_r(\la)\ \big|\ a<b,\ a\equiv d\nequiv b\ppmod e\big\},\\
N_{r,k}(\xi) &= \big\{a,b\in\mathfrak{B}_r(\xi)\ \big|\ a<b,\ a\equiv d\nequiv b\ppmod e\big\}.
\end{align*}

Suppose first that $b>d+ce$.  Then we have $\ux k(\la)-\ux k(\xi)=1$ by Lemma \ref{altux}.  On the other hand, both $b$ and $b-e$ are greater than all of $d,d+e,\dots,d+(c-1)e$, so we have
\begin{align*}
N_{r,k}(\xi) = N_{r,k}(\la)\cup&\big\{(d,b-e),(d+e,b-e),\dots,(d+(c-1)e,b-e)\big\}\\
\setminus&\big\{(d,b),(d+e,b),\dots,(d+(c-1)e,b)\big\}\end{align*}
and $n_{r,k}(\la)=n_{r,k}(\xi)$.

Alternatively, suppose $b<d+ce$.  Then $\ux k(\la)=\ux k(\xi)$ by Lemma \ref{altux}.  Let $a$ be the largest integer congruent to $d$ modulo $e$ which is less than $b$; that is, $a=b-((b-d)\mmod e)$.  Since $b<d+ce$, we have $a\in \mathfrak{B}_r(\la)$ and $a\in \mathfrak{B}_r(\xi)$, so
\begin{align*}
N_{r,k}(\xi) = N_{r,k}(\la)\cup&\big\{(d,b-e),(d+e,b-e),\dots,(a-e,b-e)\big\}\\
\setminus&\big\{(d,b),(d+e,b),\dots,(a,b)\big\}\end{align*}
and $n_{r,k}(\la)=n_{r,k}(\xi)+1$.
\end{pf}

This implies the following.

\begin{cory}\label{uxn}
Suppose $\la$ and $\mu$ are $k$-empty partitions with the same $e$-core and the same $e$-weight.  Then $\ux k(\la)=\ux k(\mu)$ if and only if $n_{r,k}(\la)=n_{r,k}(\mu)$.
\end{cory}

\sect{Canonical bases and $q$-decomposition numbers}\label{canbas}

In this section, we prove Theorem \ref{main}.  In order to do this, we must give a detailed description of how to compute the $q$-decomposition numbers.  The method we use here is via a direct computation of the bar involution.

Our set-up is largely based on \cite{lt}, where the Fock space is realised as the space of \emph{semi-infinite wedges} (of a fixed charge) modulo ordering relations; this realisation is due to Kashiwara, Miwa and Stern \cite{kms}.  Our treatment actually uses finite wedges; for each $r$, we define the truncated Fock space $\calf^e_r$ to be the span of wedges of length $r$ modulo the ordering relations, and we define a bar involution on $\calf^e_r$.  Given a partition $\mu$, there is a corresponding basis element $\bra\mu$ in $\calf^e_r$ for any $r\gs\mu'_1$; moreover, the image of $\bra\mu$ under the bar involution, when written as a linear combination of basis elements $\bra\la$, is independent of $r$, provided $r$ is sufficiently large (in fact, $r\gs|\mu|$ is sufficient; this stability result is implicit in the description in \cite{lt} of the bar involution).  So for a given partition $\mu$ we can define $\ol{\bra\mu}$ in the Fock space $\calf$ by taking a value of $r$ which is large relative to $\mu$.  This defines the bar involution on the whole of the Fock space, and hence the canonical basis and the $q$-decomposition numbers.

\subs{The truncated Fock space and the bar involution}

Fix a positive integer $r$, and define an \emph{$r$-wedge} to be a symbol of the form
\[\uu {i_1}\wed\uu {i_2}\wed\dots\wed\uu {i_r},\]
where $i_1,\dots,i_r$ are non-negative integers.  The \emph{$r$-wedge space} is the $\bbq(q)$-vector space with the set of all $r$-wedges as a basis.  We say that the $r$-wedge $\uu{i_1}\wed\dots\wed\uu{i_r}$ is \emph{ordered} if $i_1>\dots>i_r$.  Given $e\gs2$, we impose commutation relations on the $r$-wedge space, depending on our fixed integer $e\gs2$,  as follows.  First suppose that $r=2$ and $l,m$ are non-negative integers with $l\ls m$.  If $l\equiv m\ppmod e$, then we set
\[\uu l\wed\uu m = -\uu m\wed\uu l.\]
If $l\nequiv m\ppmod e$, then we define $i=(m-l)\mmod e$, and set
\begin{align*}
\uu l\wed\uu m = -q^{-1}\uu m\wed\uu l+(q^{-2}-1)\Big(&\uu{m-i}\wed\uu{l+i}\\
&-q^{-1}\uu{m-e}\wed\uu{l+e}\\
&+q^{-2}\uu{m-e-i}\wed\uu{l+e+i}\\
&-q^{-3}\uu{m-2e}\wed\uu{l+2e}\\
&+ \dots\Big),
\end{align*}
where the summation on the right continues as long as the terms are ordered.  For $r>2$, we impose the above commutation relations in every adjacent pair of positions.  The \emph{truncated Fock space} $\mathcal{F}^e_r$ is defined to be the $r$-wedge space modulo the commutation relations.

Now suppose $\mu$ is a partition, and $r\gs\mu'_1$.  Write $\mathfrak{B}_r(\mu)=\{\beta_1,\dots,\beta_r\}$ with $\beta_1>\dots>\beta_r$, and define
\[\bra\mu = \uu{\beta_1}\wed\dots\wed\uu{\beta_r}.\]
Clearly, any ordered $r$-wedge has the form $\bra\la$ for some partition $\la$ with $\la'_1\ls r$.  Any unordered $r$-wedge can be uniquely written as a linear combination of ordered $r$-wedges, so the elements $\bra\la$ with $\la'_1\ls r$ form a basis for $\mathcal{F}^e_r$, which we call the \emph{standard basis}.

Now we can define the bar involution.  Given a partition $\mu$ with $\mu'_1\ls r$, let $\bra\mu=\uu{\beta_1}\wed\dots\wed\uu{\beta_r}$ as above, and let $\arb\mu$ be the reversed wedge $\uu{\beta_r}\wed\dots\wed\uu{\beta_1}$.  Write $\arb\mu$ as a linear combination of ordered wedges using the commutation relations:
\[\arb\mu = \sum_\la b_{\la\mu}(q)\bra\la.\]
The coefficient $b_{\mu\mu}(q)$ is easy to compute (an expression is given in \cite[\S3]{lt}) and in particular is non-zero.  So we can normalise by defining $a_{\la\mu}(q) = \left.b_{\la\mu}(q)\right/b_{\mu\mu}(q)$ for all $\la,\mu$, and setting
\[\ol{\bra\mu} = \sum_\la a^e_{\la\mu}(q)\bra\la.\]

\begin{eg}
Take $e=3$, $r=4$ and $\mu=(4)$.  Then $\arb\mu = \uu 0\wed\uu1\wed\uu2\wed\uu7$, and applying the commutation relations we find that this equals
\[q^{-5}\uu7\wed\uu2\wed\uu1\wed\uu0\ +\ (q^{-4}-q^{-6})\uu5\wed\uu4\wed\uu1\wed\uu0\ +\ (q^{-7}+q^{-5})\uu4\wed\uu3\wed\uu2\wed\uu1.\]
Multiplying by $q^5$, we obtain
\[\ol{\bra{(4)}}\ =\ \bra{(4)}\ +\ (q-q^{-1})\bra{(2^2)}\ +\ (q^{-2}-1)\bra{(1^4)}.\]
\end{eg}

As mentioned above, the coefficients $a^e_{\la\mu}(q)$ are independent of the choice of $r$, provided $r$ is sufficiently large.  So we may define the bar involution on the full Fock space $\calf$: for each partition $\mu$ we define $\ol{\bra\mu}$ by computing the coefficients $a_{\la\mu}(q)$ in $\calf^e_r$ for sufficiently large $r$; then we extend semi-linearly to the whole of $\calf$, i.e.\ given coefficients $c_\mu(q)\in\bbq(q)$ we set
\[\ol{\sum_\mu c_\mu(q)\bra\mu} = \sum_\mu c_\mu(q^{-1})\ol{\bra\mu}.\]

As mentioned in \S\ref{canoni} (and as we shall shortly prove), the coefficient $a^e_{\la\mu}(q)$ is zero unless $\mu\dom\la$, and by construction the coefficient $a^e_{\mu\mu}(q)$ equals $1$.  So the canonical basis of $\calf$ and the $q$-decomposition numbers $d^e_{\la\mu}(q)$ may be defined as in Theorem \ref{bareu}, and there is a straightforward algorithm to compute them.  From this construction of the canonical basis, it follows that $d^e_{\la\mu}(q)=0$ unless $\mu\dom\la$, and that $d^e_{\la\mu}(q)$ depends only on the coefficients $a^e_{\rho\xi}(q)$ for $\mu\dom\xi\dom\rho\dom\la$; the proof of Theorem \ref{main} essentially rests on these statements.

Now we prove the promised results concerning the dominance order.  Recall from \S\ref{domsec} the partial order $\cgs$ and the definition of the $e$-extension of a multiset.

\begin{lemma}\label{strdom}
Suppose $\uu{i_1}\wed\dots\wed\uu{i_r}$ is any $r$-wedge, and write it as a linear combination of ordered wedges:
\[\uu{i_1}\wed\dots\wed\uu{i_r} = \sum_{j_1>\dots>j_r}c_{j_1\dots j_r}(q)\uu{j_1}\wed\dots\wed\uu{j_r}.\]
Then for any $j_1>\dots>j_r$ with $c_{j_1\dots j_r}\neq0$, the following statements hold.
\begin{enumerate}
\item
$j_1\ls\max\{i_1,\dots,i_r\}$ and $j_r\gs \min\{i_1,\dots,i_r\}$.
\item
The multisets
\[\{i_1\mmod e,\dots,i_r\mmod e\},\qquad \{j_1\mmod e,\dots,j_r\mmod e\}\]
are equal.
%$\big|\{k\mid i_k\equiv d\ppmod e\}\big|=\big|\{k\mid j_k\equiv d\ppmod e\}\big|$ for any $d\in\{0,\dots,e-1\}$;
\item
$X^e(i_1,\dots,i_r)\cgs X^e(j_1,\dots,j_r)$.
\end{enumerate}
\end{lemma}

\begin{pf}
Define $\inc(i_1,\dots,i_r) = \sum_{1\ls k<l\ls r}(i_k-i_l-1)^2$, and proceed by induction on $\inc(i_1,\dots,i_r)$.  If $\uu{i_1}\wed\dots\wed\uu{i_r}$ is ordered then the lemma is trivial, so suppose otherwise and choose $k$ such that $i_k\ls i_{k+1}$.  If $i_k=i_{k+1}$, then $\uu{i_1}\wed\dots\wed\uu{i_r}=0$ and again there is nothing to prove; so we can assume $i_k<i_{k+1}$.  Applying the commutation relations in positions $k,k+1$, we can write $\uu{i_1}\wed\dots\wed\uu{i_r}$ as a linear combination of wedges of the form
\[\uu{i_1}\wed\dots\wed\uu{i_{k-1}}\wed\uu{j_k}\wed\uu{j_{k+1}}\wed\uu{i_{k+2}}\wed\dots\wed\uu{i_r}\]
with $j_k>j_{k+1}$.
To prove the lemma, it suffices to show that for any such $j_k,j_{k+1}$ we have:
\begin{enumerate}
\setcounter{enumi}{-1}
\item\label{pro0}
$\inc(i_1,\dots,i_{k-1},j_k,j_{k+1},i_{k+2},\dots,i_r)<\inc(i_1,\dots,i_r)$;
\item\label{pro1}
 $\max\{i_1,\dots,i_{k-1},j_k,j_{k+1},i_{k+2},\dots,i_r\}\ls \max\{i_1,\dots,i_r\}$ and\\
$\min\{i_1,\dots,i_{k-1},j_k,j_{k+1},i_{k+2},\dots,i_r\}\gs \min\{i_1,\dots,i_r\}$;
\item\label{pro2}
$j_k,j_{k+1}$ are congruent to $i_k,i_{k+1}$ in some order, modulo $e$.
\item\label{pro3}
$X^e(i_1,\dots,i_r)\cgs X^e(i_1,\dots,i_{k-1},j_k,j_{k+1},i_{k+2},\dots,i_r)$.
\end{enumerate}
(\ref{pro0}) is a simple exercise in inequalities, using the facts
\[i_k+i_{k+1}=j_k+j_{k+1}\qquad\text{and}\qquad i_{k+1}\gs j_k>j_{k+1}\gs i_k\tag*{($\dagger$)}\]
which are immediate from the commutation relations.  (\ref{pro1}) also follows from ($\dagger$), and (\ref{pro2}) is inherent in the commutation relations.  So we are left with (\ref{pro3}).  We suppose $j_k\equiv i_k\ppmod e$; the case $j_k\equiv i_{k+1}$ is similar.  ($\dagger$) and (\ref{pro2}) imply that the  multiset $X^e(i_1,\dots,i_{k-1},j_k,j_{k+1},i_{k+2},\dots,i_r)$ may be obtained from $X^e(i_1,\dots,i_r)$ by adding a copy of each of the integers $i_k+e,i_k+2e,\dots,j_k$, and removing a copy of each of $j_{k+1}+e,j_{k+1}+2e,\dots,i_{k+1}$.  We have
\begin{align*}
i_k+e &\ls j_{k+1}+e,\\
i_k+2e&\ls j_{k+1}+2e,\\
&\ \ \vdots\\
j_k&\ls i_{k+1},
\end{align*}
and (\ref{pro3}) follows.
\end{pf}

Given our definition of the dominance order, part (3) of Lemma \ref{strdom} immediately gives the following.

\begin{cory}
Suppose $\la$ and $\mu$ are partitions.  Then the coefficient $a^e_{\la\mu}(q)$ equals $0$ unless $\mu\dom\la$.
\end{cory}

Now the following result is deduced in exactly the same way as for the standard dominance order \cite[\S4]{lt}.

\begin{propn}\label{qdecdom}
Suppose $\la$ and $\mu$ are partitions.  Then the $q$-decomposition number $d^{e}_{\la\mu}(q)$ equals $0$ unless $\mu\dom\la$.
\end{propn}

\subs{Runner removal and the bar involution}

The proof of Theorem \ref{main} will reduce to the following, which is the corresponding statement for the coefficients $a^e_{\la\mu}(q)$.

\begin{propn}\label{mainbar}
Suppose $e\gs3$, $\la$ and $\mu$ are partitions and $k\in\{0,\dots,e-1\}$.  If $\la$ and $\mu$ are $k$-empty and $\ux k(\la)=\ux k(\mu)$, then $a^{e}_{\la\mu}(q)=a^{e-1}_{\la^{-k}\mu^{-k}}(q)$.
\end{propn}

The proof of Proposition \ref{mainbar} amounts to comparing the computations of $a^{e}_{\la\mu}(q)$ and $a^{e-1}_{\la^{-k}\mu^{-k}}(q)$ using the commutation relations.  We begin by proving all the intermediate results we need concerning the commutation relations.

\begin{lemma}\label{simplewed}
%checked 27/xi
Suppose $l<m-2e$, and put $i=m-l\mmod e$.  If $i\neq0$, then:
\begin{enumerate}
\item
$\uu m\wed\uu l - \uu{m-e}\wed\uu{l+e} = -q\uu l\wed\uu m+q^{-1}\uu{l+e}\wed\uu{m-e}+(q^{-1}-q)\uu{m-i}\wed\uu{l+i}$;
\item
$\uu l\wed\uu m-\uu{l+e}\wed\uu{m-e} = -q^{-1}\uu m\wed\uu l+q\uu{m-e}\wed\uu{l+e}+(q-q^{-1})\uu{l+i}\wed\uu{m-i}$.
\end{enumerate}
\end{lemma}

\begin{pf}
Both statements are straightforward consequences of the commutation relations.
\end{pf}

\begin{lemma}\label{wedgemod}
%checked 27/xi
Suppose $l<m$ and put $i=m-l\mmod e$.  If $i\neq0$, then
\begin{enumerate}
\item\label{wedgemod1}
$\uu m\wed\uu l = -q\uu l\wed\uu m+(q^{-1}-q)\Big(\uu{l+e}\wed\uu{m-e}+\uu{l+2e}\wed\uu{m-2e}+\dots+\uu{m-i}\wed\uu{l+i}\Big)$;
\item\label{wedgemod2}
$\uu l\wed\uu m = -q^{-1}\uu m\wed\uu l+(q-q^{-1})\Big(\uu{m-e}\wed\uu{l+e}+\uu{m-2e}\wed\uu{l+2e}+\dots+\uu{l+i}\wed\uu{m-i}\Big)$.
\end{enumerate}
\end{lemma}

\begin{pf}
Both statements are easily proved by induction on $m-l$; the cases $m-l=i$ and $m-l=e+i$ follow easily from the  commutation relations, and the inductive step from Lemma \ref{simplewed}.
\end{pf}

Now we fix $d\in\{0,\dots,e-1\}$.

\begin{lemma}\label{oned}
Suppose $h_1,\dots,h_s,i\in\zo$, with $h_1,\dots,h_s\nequiv d\equiv i\ppmod{e}$.  Then the wedge $\uu{h_1}\wed\dots\wed\uu{h_s}\wed\uu i$ can be expressed as a linear combination of wedges $\uu j\wed\uu{k_1}\wed\dots\wed\uu{k_s}$ in which $j\equiv d\nequiv k_1,\dots,k_s\ppmod{e}$, and $j,k_1,\dots,k_s\ls\max\{h_1,\dots,h_s,i\}$.

Furthermore, if $\max\{h_1,\dots,h_s,i\}=i$, then we may construct such a linear combination in such a way that the only wedge of the form $\uu i\wed\uu{k_1}\wed\dots\wed\uu{k_s}$ occurring with non-zero coefficient is $\uu{i}\wed\uu{h_1}\wed\dots\wed\uu{h_s}$.
\end{lemma}

\begin{pf}
We proceed by induction on $s$.  The case $s=0$ is trivial, so assume $s\gs1$.  Using Lemma \ref{wedgemod}, we may express $\uu{h_s}\wed\uu{i}$ as a linear combination of wedges $\uu{i'}\wed\uu{k}$, with $i'\equiv d\nequiv k\ppmod{e}$; we use part (\ref{wedgemod1}) of that lemma if $h_s>i$, or part (\ref{wedgemod2}) if $h_s<i$.  Either way, we see that for every such wedge we have $i',k\ls\max\{h_s,i\}$; moreover, if $h_s<i$, then the wedge $\uu i\wed\uu {h_s}$ occurs with non-zero coefficient (and no other wedge of the form $\uu i\wed\uu k$ occurs).  Now for each such pair $(i',k)$, we apply the inductive hypothesis to the $s$-wedge $\uu{h_1}\wed\dots\wed\uu{h_{s-1}}\wed\uu{i'}$.
\end{pf}

\begin{cory}\label{manyd}
Suppose we are given $h_1,\dots,h_s,i_1,\dots,i_c\in\zo$ and $0\ls m_1\ls \dots\ls m_c\ls s$ such that:
\begin{itemize}
\item
$i_1,\dots,i_c\equiv d\pmod{e}$ and $i_1<\dots<i_c$;
\item
$h_1,\dots,h_s\nequiv d\pmod{e}$;
\item
$i_x>h_1,\dots,h_{m_x}$ for each $x\in\{1,\dots,c\}$.
\end{itemize}
Then the $(s+c)$-wedge
\[\Big(\uu{h_1}\wed\dots\wed\uu{h_{m_1}}\Big)\wed\uu{i_1}\wed\Big(\uu{h_{m_1+1}}\wed\dots\wed\uu{h_{m_2}}\Big)\wed\uu{i_2}\wed\quad\dots\quad\wed\uu{i_c}\wed\Big(\uu{h_{m_c+1}}\wed\dots\wed\uu{h_s}\Big)\]
may be written as a linear combination of wedges of the form
\[\uu{j_1}\wed\dots\wed\uu{j_c}\wed\uu{k_1}\wed\dots\wed\uu{k_s},\]
where
\begin{itemize}
\item
$j_1,\dots,j_c\equiv d\pmod{e}$;
\item
$k_1,\dots,k_s\nequiv d\pmod{e}$;
\item
$j_x\ls i_x$ for each $x\in\{1,\dots,c\}$.
\end{itemize}
Furthermore, this may be done in such a way that the only wedge of the form $\uu{i_1}\wed\dots\wed\uu{i_c}\wed\uu{k_1}\wed\dots\wed\uu{k_s}$ occurring with non-zero coefficient is
\[\uu{i_1}\wed\dots\wed\uu{i_c}\wed\uu{h_1}\wed\dots\wed\uu{h_s}.\]
\end{cory}

\begin{pf}
We use induction on $c$, with the case $c=0$ being trivial.  Assuming $c\gs1$, we apply Lemma \ref{oned} to the $(m_1+1)$-wedge $\uu{h_1}\wed\dots\wed\uu{h_{m_1}}\wed\uu{i_1}$.  This yields a linear combination of wedges of the form $\uu{j_1}\wed \uu{k_1}\wed\dots\wed\uu{k_{m_1}}$, with $j_1,k_1,\dots,k_{m_1}\ls i_1$ and $j_1\equiv d\ppmod{e}$; furthermore, the only such wedge occurring with non-zero coefficient in which $j_1=i_1$ is the wedge $\uu{i_1}\wed\uu{h_1}\wed\dots\wed\uu{h_{m_1}}$.

Given a wedge $\uu{j_1}\wed \uu{k_1}\wed\dots\wed\uu{k_{m_1}}$ occurring in this linear combination, we have $k_1,\dots,k_{m_1}\ls i_1<i_2$, so we can apply the inductive hypothesis to the $(s+c-1)$-wedge
\[\Big(\uu{k_1}\wed\dots\wed\uu{k_{m_1}}\wed\uu{h_{m_1+1}}\wed\dots\wed\uu{h_{m_2}}\Big)\wed\uu{i_2}\wed\quad\dots\quad\wed\uu{i_c}\wed\Big(\uu{h_{m_c+1}}\wed\dots\wed\uu{h_s}\Big),\]
which gives the result.
\end{pf}

\begin{cory}\label{shoveleft}
Suppose we are given $h_1,\dots,h_s\in\zo$ and $0\ls m_1\ls \dots\ls m_c\ls s$ such that
\begin{itemize}
\item
$h_1,\dots,h_s\nequiv d\pmod{e}$;
\item
$h_1,\dots,h_{m_x}<d+(x-1)e$ for each $x\in\{1,\dots,c\}$.
\end{itemize}
Let $w$ be the $(s+c)$-wedge
\[\Big(\uu{h_1}\wed\dots\wed\uu{h_{m_1}}\Big)\wed\uu{d}\wed\Big(\uu{h_{m_1+1}}\wed\dots\wed\uu{h_{m_2}}\Big)\wed\uu{d+e}\wed\quad\dots\quad\wed\uu{d+(c-1)e}\wed\Big(\uu{h_{m_c+1}}\wed\dots\wed\uu{h_s}\Big).\]
Then $w$ equals a non-zero multiple of the wedge
\[\Big(\uu d\wed\uu{d+e}\wed\dots\wed\uu{d+(c-1)e}\Big)\wed\Big(\uu{h_1}\wed\dots\wed\uu{h_s}\Big).\]
\end{cory}

\begin{pf}
Apply Corollary \ref{manyd} to $w$, and suppose $\uu{j_1}\wed\dots\wed\uu{j_c}\wed\uu{k_1}\wed\dots\wed\uu{k_s}$ is one of the resulting wedges.  Since $j_1,\dots,j_c$ are all congruent modulo $e$, we may re-write $\uu{j_1}\wed\dots\wed\uu{j_c}$ as $\pm\uu{j_{\pi(1)}}\wed\dots\wed\uu{j_{\pi(c)}}$, where $\pi\in\sss c$ is such that $j_{\pi(1)}\gs\dots\gs j_{\pi(c)}$.  If any two of $j_1,\dots,j_c$ are equal, then the latter wedge will equal zero.  So we may discard any terms in which $j_1,\dots,j_c$ are not pairwise distinct.  But recall that we have $j_1,\dots,j_c\equiv d\ppmod{e}$, and $j_x\ls d+(x-1)e$ for each $x$.  The only way such $j_1,\dots,j_c$ can be pairwise distinct is if $j_x=d+(x-1)e$ for each $x$.  Now the last statement of Corollary \ref{manyd} gives the result.
\end{pf}

\begin{eg}
Suppose $e=3$ and $\mu=(7,3,1)$.  Taking $r=6$, we get $\mathfrak{B}_r(\mu) = \{12,7,4,2,1,0\}$, so that
\[\arb{\mu} = \uu 0\wed\uu1\wed\uu2\wed\uu4\wed\uu7\wed\uu{12}.\]
Taking $d=1$ and applying Lemma \ref{wedgemod} repeatedly to move terms congruent to $1$ modulo $3$ to the left, we find that
\begin{alignat*}3
\arb{\mu} =&{}& &-q^{-5}&&\uu1\wed\uu4\wed\uu7\wed\uu0\wed\uu2\wed\uu{12}\\
&&+q^{-4}&(q-q^{-1})&&\uu1\wed\uu4\wed\uu4\wed\uu3\wed\uu2\wed\uu{12}\\
&&+q^{-4}&(q-q^{-1})&&\uu1\wed\uu4\wed\uu1\wed\uu6\wed\uu2\wed\uu{12}\\
&&+q^{-4}&(q-q^{-1})&&\uu1\wed\uu1\wed\uu7\wed\uu3\wed\uu2\wed\uu{12}\\
&&-q^{-3}&(q-q^{-1})^2&&\uu1\wed\uu1\wed\uu4\wed\uu6\wed\uu2\wed\uu{12}\\
&&+q^{-4}&(q-q^{-1})&&\uu1\wed\uu4\wed\uu4\wed\uu0\wed\uu5\wed\uu{12}\\
&&-q^{-3}&(q-q^{-1})^2&&\uu1\wed\uu4\wed\uu1\wed\uu3\wed\uu5\wed\uu{12}\\
&&-q^{-3}&(q-q^{-1})^2&&\uu1\wed\uu1\wed\uu4\wed\uu3\wed\uu5\wed\uu{12}.
\end{alignat*}
When we apply the commutation relations in the first three positions, all terms apart from the first vanish, so that $\arb{\mu}$ equals a non-zero multiple of $\uu1\wed\uu4\wed\uu7\wed\uu0\wed\uu2\wed\uu{12}$.
\end{eg}

Now we need to compare the Fock spaces $\mathcal{F}^e_s$ and $\mathcal{F}^{e-1}_s$.  To avoid ambiguity, we write a wedge in the latter Fock space as
\[\uu{i_1}\wede\dots\wede\uu{i_s};\]
so wedges written in this way are subject to the commutation relations modulo $e-1$, while wedges written using the symbol $\wed$ are subject to the commutation relations modulo $e$.  Now, recalling the function $\phi_d$ from \S\ref{31}, we have the following.

\begin{lemma}\label{fstraight}
Suppose we have $i_{kl}\in\zo$ for $1\ls k\ls t$ and $1\ls l\ls s$, with $i_{kl}\nequiv d\ppmod e$, and suppose $b_1(q),\dots,b_t(q)\in\bbq(q)$ are such that
\[\sum_{k=1}^t b_k(q)\Big(\uu{i_{k1}}\wed\dots\wed\uu{i_{ks}}\Big)=0.\]
Then
\[\sum_{k=1}^t b_k(q)\Big(\uu{\phi_d(i_{k1})}\wede\dots\wede\uu{\phi_d(i_{ks})}\Big)=0.\]

\end{lemma}

\begin{pf}
This comes directly from a comparison of the commutation relations for $\uu l\wed\uu m$ and for $\uu{\phi_d(l)}\wede\uu{\phi_d(m)}$ when $l\ls m$.  We leave this for the reader to check.
\end{pf}

\begin{eg}
Suppose $e=3$.  Then the commutation relations give
\begin{align*}
\uu0\wed\uu2\wed\uu{12} = &-q^{-2}\uu{12}\wed\uu2\wed\uu0+(q^{-3}-q^{-1})\uu{11}\wed\uu3\wed\uu0\\
&+(q^{-2}-q^{-4})\uu9\wed\uu5\wed\uu0+(q^{-5}-q^{-3})\uu8\wed\uu6\wed\uu0,
\end{align*}
while the commutation relations modulo $2$ give
\begin{align*}
\uu0\wede\uu1\wede\uu8 = &-q^{-2}\uu8\wede\uu1\wede\uu0+(q^{-3}-q^{-1})\uu7\wede\uu2\wede\uu0\\
&+(q^{-2}-q^{-4})\uu6\wede\uu3\wede\uu0+(q^{-5}-q^{-3})\uu5\wede\uu4\wede\uu0.
\end{align*}
\end{eg}

We need one more ingredient before we prove our main result; this needs some preparation.

\begin{lemma}\label{shovr0}
Suppose $j\in\zo$ and $w=\uu{i_1}\wed\dots\wed\uu{i_u}$ is a $u$-wedge having the following property: there is a unique $x\in\{1,\dots,u\}$ such that $i_x\equiv d\ppmod e$, and for this value of $x$ we have $i_1,\dots,i_x>j$.  Then, when we write $w$ as a linear combination of ordered wedges using the commutation relations, every ordered wedge that occurs with non-zero coefficient contains exactly one term $\uu i$ with $i\equiv d\ppmod e$, and this value of $i$ satisfies $i>j$.
\end{lemma}

\begin{pf}
If $i_1\gs\dots\gs i_u$ then there is nothing to prove, so we suppose $i_y<i_{y+1}$ for some $y$, and apply the commutation relations in positions $y,y+1$.  This gives an expression for $\uu{i_1}\wed\dots\wed\uu{i_u}$ as a linear combination of wedges of the form $\uu{i_1}\wed\dots\wed\uu{i_{y-1}}\wed\uu{l_y}\wed\uu{l_{y+1}}\wed\uu{i_{y+2}}\wed\dots\wed\uu{i_u}$.  Defining $\inc(i_1,\dots,i_u)$ as in the proof of Lemma \ref{strdom}, it suffices to show that $\inc(i_1,\dots,i_{y-1},l_y,l_{y+1},i_{y+2},\dots,i_u)<\inc(i_1,\dots,i_u)$ and that the hypotheses of the lemma hold with $i_1,\dots,i_u$ replaced by $i_1,\dots,i_{y-1},l_y,l_{y+1},i_{y+2},\dots,i_u$.  The first fact is a simple exercise as before, and the second fact is easy to check from the commutation relations.
\end{pf}

\begin{lemma}\label{shovr1}
Suppose $j,k_1,\dots,k_s\in\zo$ with $j\equiv d\nequiv k_1,\dots,k_s\ppmod{e}$ and $k_1>\dots> k_s> j$.  Then, when the $(s+1)$-wedge $w=\uu j\wed\uu{k_1}\wed\dots\wed\uu{k_s}$ is expressed as a linear combination of ordered wedges using the commutation relations, each wedge that occurs contains exactly one term $\uu{i}$ with $i\equiv d\ppmod{e}$, and this term satisfies $i\gs j$.  Moreover, the only wedge occurring that includes the term $\uu{j}$ is the wedge
\[\uu{k_1}\wed\dots\wed\uu{k_s}\wed\uu{j},\]
occurring with coefficient $(-q)^{-s}$.
\end{lemma}

\begin{pf}
We use induction on $s$, with the case $s=0$ being trivial.  Assuming $s\gs1$, we apply the commutation relations in positions $1$ and $2$.  This yields an expression
\[w = -q^{-1}w_1+\sum_{v\in V}b_v(q)v,\]
where $w_1 = \uu{k_1}\wed\uu j\wed\uu{k_2}\wed\dots\wed\uu{k_s}$, and $V$ is a set of wedges each of which satisfies the hypotheses of Lemma \ref{shovr0}.  Applying Lemma \ref{shovr0} to any $v\in V$, we get a linear combination of ordered wedges in which there is one term $\uu i$ with $i\equiv d$, and this $i$ satisfies $i>j$.  So we can neglect all wedges $v\in V$, and it suffices to show that the present lemma holds with $w_1$ in place of $w$ and $(-q)^{1-s}$ in place of $(-q)^{-s}$.

By induction, when we write
\[\uu j\wed\uu{k_2}\wed\dots\wed\uu{k_s}=\sum_{j_1>\dots>j_s}c_{j_1\dots j_s}(q)\uu{j_1}\wed\dots\wed\uu{j_s}\]
we have that:
\begin{itemize}
%\item
%each $\uu{j_1}\wed\dots\wed\uu{j_s}$ is ordered;
\item
if $c_{j_1\dots j_s}(q)\neq0$, then there is exactly one $x$ such that $j_x\equiv d\ppmod e$, and this $j_x$ is greater than or equal to $j$;
\item
if $c_{j_1\dots,j_s}(q)\neq0$ and $j_x=j$ for some $x$, then $(j_1,\dots,j_s)=(k_2,\dots,k_s,j)$ and $c_{j_1\dots j_s}(q) = (-q)^{1-s}$.
\end{itemize}
Also, by Lemma \ref{strdom}(\ref{pro1}), each $(j_1,\dots,j_s)$ with $c_{j_1\dots,j_s}(q)\neq0$ satisfies $j_1\ls k_2<k_1$.  So we see that
\[w_1=\sum_{j_1,\dots,j_s}c_{j_1\dots j_s}(q)\uu{k_1}\wed\uu{j_1}\wed\dots\wed\uu{j_s}\]
is an expression for $w_1$ as a linear combination of ordered wedges, and the result follows.
\end{pf}

\begin{lemma}\label{shovr2}
Suppose $j,k_1,\dots,k_s\in\zo$ with $j\equiv d\ppmod e$ and $k_1>\dots>k_s$.  Suppose also that for some $z\in\{1,\dots,s\}$ we have $k_z\equiv d\ppmod e$ and $k_z\gs j$.  Let $w$ denote the $(s+1)$-wedge $\uu j\wed\uu{k_1}\wed\dots\wed\uu{k_s}$.  When $w$ is written as a linear combination of ordered wedges using the commutation relations, every wedge that occurs with non-zero coefficient contains a term $\uu l$ with $l\equiv d\ppmod e$ and $l\gs j+e$.
\end{lemma}

\begin{pf}
Let $y$ be maximal such that $k_y\gs j$, and let $w'$ be the $(y+1)$-wedge $\uu j\wed\uu{k_1}\wed\dots\wed\uu{k_y}$.  Write $w'$ as a linear combination of ordered wedges:
\[w' = \sum_{l_1,\dots,l_{y+1}}b_{l_1\dots l_{y+1}}(q)\uu{l_1}\wed\dots\wed\uu{l_{y+1}}.\]
If $b_{l_1\dots l_{y+1}}(q)\neq0$, then by Lemma \ref{strdom}(\ref{pro1}) we have $l_{y+1}\gs\min\{j,k_1,\dots,k_y\}=j$, and by Lemma \ref{strdom}(\ref{pro2}) at least two of $l_1,\dots,l_{y+1}$ are congruent to $d$ modulo $e$.  Hence for some $x\ls y$ we have $l_x\equiv d\ppmod e$ and $l_x\gs j+e$.  The fact that $l_{y+1}\gs j>k_{y+1}$ implies that the $(s+1)$-wedge
\[\uu{l_1}\wed\dots\wed\uu{l_{y+1}}\wed\uu{k_{y+1}}\wed\dots\wed\uu{k_s}\]
is ordered.  So we see that
\[w = \sum_{l_1,\dots,l_{y+1}}b_{l_1\dots l_{y+1}}(q)\uu{l_1}\wed\dots\wed\uu{l_{y+1}}\wed\uu{k_{y+1}}\wed\dots\wed\uu{k_s}\]
is an expression for $w$ as a linear combination of ordered wedges with the required properties.
\end{pf}

\begin{cory}\label{shoveright}
Suppose $k_1,\dots,k_s\in\zo$, with $k_1,\dots,k_s\nequiv d\ppmod{e}$ and $k_1>\dots> k_s$, and let $0\ls n_1\ls \dots\ls  n_c\ls s$ be such that $k_{n_x}>d+(c-x)e>k_{n_x+1}$ for all $x$.  Then, when the $(s+c)$-wedge
\[\Big(\uu{d+(c-1)e}\wed\uu{d+(c-2)e}\wed\dots\wed\uu d\Big)\wed\Big(\uu{k_1}\wed\dots\wed\uu{k_s}\Big)\]
is expressed as a linear combination of ordered wedges, the only ordered wedge occurring that includes all the terms $\uu{d+(c-1)e},\uu{d+(c-2)e},\dots,\uu d$ is the wedge
\[\Big(\uu{k_1}\wed\dots\wed\uu{{k_{n_1}}}\Big)\wed\uu{d+(c-1)e}\wed\Big(\uu{k_{n_1+1}}\wed\dots\wed\uu{k_{n_2}}\Big)\wed\uu{d+(c-2)e}\wed\quad\dots\quad\wed\uu{d}\wed\Big(\uu{k_{n_c+1}}\wed\dots\wed\uu{k_s}\Big),\]
occurring with coefficient $(-q)^{-(n_1+\dots+n_s)}$.
\end{cory}

\begin{pf}
We use induction on $c$, with the case $c=0$ being trivial.  Assuming $c\gs 1$, apply the inductive hypothesis to write
\[\Big(\uu{d+(c-2)e}\wed\uu{d+(c-3)e}\wed\dots\wed\uu d\Big)\wed\Big(\uu{k_1}\wed\dots\wed\uu{k_s}\Big) = \sum_u b_u(q)u,\]
where:
\begin{itemize}
\item
each $u$ is an ordered wedge;
\item
if we set
\[w'=\Big(\uu{k_1}\wed\dots\wed\uu{k_{n_2}}\Big)\wed\uu{d+(c-2)e}\wed\Big(\uu{k_{n_2+1}}\wed\dots\wed\uu{k_{n_3}}\Big)\wed\uu{d+(c-3)e}\wed\quad\dots\quad\wed\uu{d}\wed\Big(\uu{k_{n_c+1}}\wed\dots\wed\uu{k_s}\Big),\]
then $b_{w'}(q)=(-q)^{-(n_2+\dots+n_s)}$;
\item
if $u\neq w'$ and $b_u(q)\neq0$, then $u$ does not contain all of the terms $\uu{d+(c-2)e},\uu{d+(c-3)e},\dots,\uu d$.
\end{itemize}
Now by Lemma \ref{strdom}(\ref{pro2}), every wedge $u$ with $b_u(q)\neq0$ contains exactly $c-1$ terms of the form $\uu{d+ze}$ with $z\in\zo$; if $u\neq w'$ then these terms are not $\uu{d+(c-2)e},\dots,\uu d$, so $u$ contains a term $\uu{d+ze}$ with $z\gs c-1$.  So by Lemma \ref{shovr2}, when we write $\uu{d+(c-1)e}\wed u$ as a linear combination of ordered wedges, each wedge that occurs contains a term $\uu{d+ze}$ with $z\gs c$, and therefore does not contain all the terms $\uu{d+(c-1)e},\uu{d+(c-2)e}\dots,\uu{d}$.  So we may ignore all terms $\uu{d+(c-1)e}\wed u$ with $u\neq w'$,  and we concentrate on the wedge $\uu{d+(c-1)e}\wed w'$.

Write
\[w_1=\uu{d+(c-1)e}\wed\uu{k_1}\wed\dots\wed\uu{k_{n_1}},\]
and express $w_1$ as a linear combination of ordered wedges:
\[w_1 = \sum_v c_v(q) v.\]
By Lemma \ref{shovr1}, any wedge $v$ with $c_v(q)\neq0$ contains a term $\uu{d+ze}$ with $z\gs c-1$, and if $v$ contains the term $\uu{d+(c-1)e}$ then $v=\uu{k_1}\wed\dots\wed\uu{k_{n_1}}\wed\uu{d+(c-1)e}$ and $c_v(q) = (-q)^{-n_1}$.  Moreover, if $v=\uu{v_1}\wed\dots\wed\uu{v_{n_1+1}}$ and $c_v(q)\neq0$, then by Lemma \ref{strdom}(\ref{pro1}) we have $v_{n_1+1}\gs d+(c-1)e>k_{n_1+1}$, so the wedge
\[v\wed \Big(\uu{k_{n_1+1}}\wed\dots\wed\uu{k_{n_2}}\Big)\wed\uu{d+(c-2)e}\wed\quad\dots\quad\wed\uu{d}\wed\Big(\uu{k_{n_c+1}}\wed\dots\wed\uu{k_s}\Big)\]
is ordered.  So an expression for $\uu{d+(c-1)e}\wed w'$ as a linear combination of ordered wedges is
\[\uu{d+(c-1)e}\wed w' = \sum_v c_v(q)\Bigg(v\wed \Big(\uu{k_{n_1+1}}\wed\dots\wed\uu{k_{n_2}}\Big)\wed\uu{d+(c-2)e}\wed\quad\dots\quad\wed\uu{d}\wed\Big(\uu{k_{n_c+1}}\wed\dots\wed\uu{k_s}\Big)\Bigg),\]
and the result follows.
\end{pf}

\begin{eg}
Taking $e=3$, $r=6$, $d=1$, we have
\begin{alignat*}3
\uu7\wed\uu4\wed\uu1\wed\uu{12}\wed\uu2\wed\uu0 =&{}& &q^{-4} &&\uu{12}\wed\uu7\wed\uu4\wed\uu2\wed\uu1\wed\uu0\\
&&+q^{-4}&(q-q^{-1})&&\uu{10}\wed\uu9\wed\uu4\wed\uu2\wed\uu1\wed\uu0\\
&&+q^{-3}&(q-q^{-1})&&\uu{10}\wed\uu7\wed\uu6\wed\uu2\wed\uu1\wed\uu0\\
&&-q^{-1}&(q-q^{-1})&&\uu{10}\wed\uu7\wed\uu4\wed\uu3\wed\uu2\wed\uu0
\end{alignat*}
and also
\begin{alignat*}3
\uu7\wed\uu4\wed\uu1\wed\uu{11}\wed\uu3\wed\uu0 =&{}& &q^{-4} &&\uu{11}\wed\uu7\wed\uu4\wed\uu3\wed\uu1\wed\uu0\\
&&+q^{-4}&(q-q^{-1})&&\uu{10}\wed\uu8\wed\uu4\wed\uu3\wed\uu1\wed\uu0\\
&&+q^{-3}&(q-q^{-1})&&\uu{10}\wed\uu7\wed\uu5\wed\uu3\wed\uu1\wed\uu0\\
&&+q^{-2}&(q-q^{-1})&&\uu{10}\wed\uu7\wed\uu4\wed\uu3\wed\uu2\wed\uu0.
\end{alignat*}
\end{eg}

Now we can prove Proposition \ref{mainbar}; the reader should combine the last three examples in this section to follow the proof for the case $e=3$, $\mu=(7,3,1)$, $\la=(6,3,1^2)$.

\begin{pfof}{Proposition \ref{mainbar}}
Choose a large $r$, and let $c$ and $d$ be as defined in \S\ref{absec}.  For any $k$-empty partition $\xi$ with the same $e$-core and $e$-weight as $\mu$, we know that the $r$-beta-set $\mathfrak{B}_r(\xi)$ contains the integers $d,d+e,\dots,d+(c-1)e$, together with $r-c$ integers not congruent to $d$ modulo $e$, which we write as $h_1(\xi)>\dots>h_{r-c}(\xi)$.  Then the $(r-c)$-beta-set $\mathfrak{B}_{r-c}(\xi^{-k})$ equals  $\{\phi_d(h_1(\xi)),\dots,\phi_d(h_{r-c}(\xi))\}$.

In the particular case $\xi=\mu$, we have
\begin{alignat*}1
\arb{\mu} = \Big(\uu{h_{r-c}(\mu)}\wed\dots\wed\uu{h_{m_c+1}(\mu)}\Big)&\wed\uu d\wed\Big(\uu{h_{m_c}(\mu)}\wed\dots\wed\uu{h_{m_{c-1}+1}(\mu)}\Big)\\
&\wed\uu{d+e}\wed\Big(\uu{h_{m_{c-1}}(\mu)}\wed\dots\wed\uu{h_{m_{c-2}+1}(\mu)}\Big)\\
&\ \ \vdots\\
&\wed\uu{d+(c-1)e}\wed\Big(\uu{h_{m_1}(\mu)}\wed\dots\wed\uu{h_1(\mu)}\Big),\\
\end{alignat*}
for appropriate $1\ls m_1\ls \dots\ls m_c\ls r$.  By Corollary \ref{shoveleft}, this equals a non-zero multiple of the wedge
\[\Big(\uu d\wed\uu{d+e}\wed\dots\wed\uu{d+(c-1)e}\Big)\wed\Big(\uu{h_{r-c}(\mu)}\wed\dots\wed\uu{h_1(\mu)}\Big).\]
Now we examine $\mu^{-k}$.  Using $(r-c)$-wedges, we have
\[\arb{\mu^{-k}} = \uu{\phi_d(h_{r-c}(\mu))}\wede\dots\wede\uu{\phi_d(h_1(\mu))},\]
so (from the definition of the constants $a^{e-1}_{\pi\mu^{-k}}(q)$)
\[\uu{\phi_d(h_{r-c}(\mu))}\wede\dots\wede\uu{\phi_d(h_1(\mu))} = C\sum_\pi a^{e-1}_{\pi\mu^{-k}}(q)\bra\pi\]
in the Fock space $\mathcal{F}^{e-1}_{r-c}$, for some non-zero $C$.  Each $\pi$ occurring has the same $(e-1)$-core and $(e-1)$-weight as $\mu^{-k}$, and so can be written as $\xi^{-k}$ for some $k$-empty partition $\xi$ with the same $e$-core and $e$-weight as $\mu$.  So by Lemma \ref{fstraight}, we get
\[\uu{h_{r-c}(\mu)}\wed\dots\wed\uu{h_1(\mu)} = C\sum_\xi a^{e-1}_{\xi^{-k}\mu^{-k}}(q).\uu{h_1(\xi)}\wed\dots\wed\uu{h_{r-c}(\xi)}\]
in the Fock space $\mathcal{F}^e_{r-c}$.  Combining this with the expression for $\arb\mu$ above and the fact that
\[\uu d\wed\dots\wed \uu{d+(c-1)e} = (-1)^{\binom c2}\uu{d+(c-1)e}\wed\dots\wed\uu d,\]
we find that $\arb{\mu}$ equals a non-zero multiple of
\[\sum_\xi a^{e-1}_{\xi^{-k}\mu^{-k}}(q).\Big(\uu{d+(c-1)e}\wed\dots\wed\uu d\Big)\wed\Big(\uu{h_1(\xi)}\wed\dots\wed\uu{h_{r-c}(\xi)}\Big),\]
summing over all $k$-empty partitions $\xi$ with the same $e$-core and $e$-weight as $\mu$.  Choose such a partition $\xi$, let $w_\xi$ denote the wedge
\[\Big(\uu{d+(c-1)e}\wed\dots\wed\uu d\Big)\wed\Big(\uu{h_1(\xi)}\wed\dots\wed\uu{h_{r-c}(\xi)}\Big),\]
and let $0\ls n_1\ls\dots\ls n_c\ls r-c$ be such that $h_{n_x}(\xi)>d+(c-x)e>h_{n_x+1}(\xi)$ for each $x$.  Note that $n_1+\dots+n_c$ is the integer $n_{r,k}(\xi)$ defined in \S\ref{analter}.  By Corollary \ref{shoveright}, $w_\xi$ equals $(-q)^{-(n_1+\dots+n_c)}$ times
\begin{align*}
\Big(\uu{h_1(\xi)}\wed\dots\wed\uu{{h_{n_1}(\xi)}}\Big)&\wed\uu{d+(c-1)e}\wed\Big(\uu{h_{n_1+1}(\xi)}\wed\dots\wed\uu{h_{n_2}(\xi)}\Big)\\
&\wed\uu{d+(c-2)e}\wed\Big(\uu{h_{n_2+1}(\xi)}\wed\dots\wed\uu{h_{n_3}(\xi)}\Big)\\
&\ \ \vdots\\
&\wed\uu{d}\wed\Big(\uu{h_{n_c+1}(\xi)}\wed\dots\wed\uu{h_{r-c}(\xi)}\Big)\\&=\bra\xi
\end{align*}
plus a linear combination of other wedges, none of which includes all the terms $\uu{d+(c-1)e},\dots,\uu d$.  Summing over $\xi$, we see that
\[\arb\mu= D\left(\sum_{\xi\in \mathcal{N}} a^{e-1}_{\xi^{-k}\mu^{-k}}(q)(-q)^{-n_{r,k}(\xi)}\bra{\xi}\right)+\sum_{\rho\in \partn\setminus \mathcal{N}}f_\rho(q)\bra\rho\]
where $D$ is a non-zero constant, $\mathcal{N}$ is the set of $k$-empty partitions with the same $e$-core and $e$-weight as $\mu$, and $f_\rho(q)\in\bbq(q)$ for $\rho\in \partn\setminus \mathcal{N}$.  Normalising, we see that
\[a^{e}_{\la\mu}(q) = \frac{Da^{e-1}_{\la^{-k}\mu^{-k}}(q)(-q)^{-n_{r,k}(\la)}}{Da^{e-1}_{\mu^{-k}\mu^{-k}}(q)(-q)^{-n_{r,k}(\mu)}}= (-q)^{n_{r,k}(\mu)-n_{r,k}(\la)}a^{e-1}_{\la^{-k}\mu^{-k}}(q).\]
But by Corollary \ref{uxn}, $\ux k(\la)=\ux k(\mu)$ implies that $n_{r,k}(\la)=n_{r,k}(\mu)$, and the proposition is proved.
\end{pfof}

\subs{Proof of Theorem \ref{main}}

Now we can complete the proof of Theorem \ref{main}.  We fix $\mu$, and proceed by induction on $\la$ with respect to the dominance order.  If $\la=\mu$ then the result is immediate, so assume $\la\neq\mu$.  Comparing coefficients of $\bra\la$ in the expression $G(\mu)=\ol{G(\mu)}$, we find that
\begin{align*}
d^{e}_{\la\mu}(q) - d^{e}_{\la\mu}(q^{-1}) &= \sum_{\xi\neq\la}d^{e}_{\xi\mu}(q^{-1})a^{e}_{\la\xi}(q).\\
\intertext{Similarly, we have}
d^{e-1}_{\la^{-k}\mu^{-k}}(q) - d^{e-1}_{\la^{-k}\mu^{-k}}(q^{-1}) &= \sum_{\pi\neq\la^{-k}}d^{e-1}_{\pi\mu^{-k}}(q^{-1})a^{e-1}_{\la^{-k}\pi}(q).
\end{align*}
Now $d^{e}_{\xi\mu}(q^{-1})=0$ unless $\mu\dom_e\xi$, while $a^{e}_{\la\xi}(q)=0$ unless $\xi\dom_e\la$, so we may restrict the range of summation in the first equation above to those $\xi$ such that $\mu\dom_e\xi\doms_e\la$.  Similarly, we may restrict the range of summation in the second equation to $\mu^{-k}\dom_{e-1}\pi\doms_{e-1}\la^{-k}$.  Now since $\ux k(\la)=\ux k(\mu)$, the set of $\pi$ with $\mu^{-k}\dom_{e-1}\pi\doms_{e-1}\la^{-k}$ is precisely the set of $\xi^{-k}$ for $\xi\in\partn$ with $\mu\dom_e\xi\doms_e\la$ (Lemma \ref{sandwich}).  Moreover, we know that for any such $\xi$ we have $\ux k(\xi)=\ux k(\mu)$, so we have
\[a^{e}_{\la\xi}(q) = a^{e-1}_{\la^{-k}\xi^{-k}}(q)\]
by Proposition \ref{mainbar}, and
\[d^{e}_{\xi\mu}(q^{-1}) = d^{e-1}_{\xi^{-k}\mu^{-k}}(q^{-1})\]
by induction.  We deduce that
\[d^{e}_{\la\mu}(q) - d^{e}_{\la\mu}(q^{-1}) = d^{e-1}_{\la^{-k}\mu^{-k}}(q) - d^{e-1}_{\la^{-k}\mu^{-k}}(q^{-1}),\]
and since $d^{e}_{\la\mu}(q),d^{e-1}_{\la^{-k}\mu^{-k}}(q)$ are polynomials divisible by $q$, the result follows.

\sect{The Mullineux map}\label{mullsec}

In this section, we examine the Mullineux map in detail and prove Theorem \ref{mainmull}.%; in particular, we look at the relationship between the Mullineux map and the function $\ux k$, and prove the following theorem.

%We now give our second main theorem, which concerns the Mullineux map.  Its significance may not be immediately obvious, but in Section \ref{conseq} we shall see that when combined with Theorem \ref{main} it is particularly powerful.  Let $m$ denote the Mullineux map.

%\begin{thm}\label{mainmull}
%Suppose $\mu$ is an $e$-regular partition.  If $\mu$ and $m(\mu)'$ are both $k$-empty, then $\ux k(\mu)= \ux k(m(\mu)')$.
%\end{thm}

%The significance of this result may not be immediately apparent, but in Section \ref{conseq} we shall see that in conjunction with Theorem \ref{main} it provides powerful information about the computation of $q$-decomposition numbers.

\subs{Definition of the Mullineux map}

The description of the Mullineux map that we use is based on the abacus, and largely taken from \cite{gv}.

\begin{defn}\label{mulldef}
Suppose $\mu$ is an $e$-regular partition, and take an abacus display for $\mu$ with $r$ beads, for some $r\gs\mu'_1$.  Let $\beta,\gamma$ be the positions of the last bead and the first empty space on the abacus, respectively; so $\beta$ is the beta-number $\beta_1=\mu_1+r-1$, while $\gamma$ equals $r-\mu'_1$.

Assuming $\mu\neq\varnothing$, there is a unique sequence $b_1>c_1>\dots>b_t>c_t$ of non-negative integers satisfying the following conditions.
\begin{enumerate}
\item
For each $1\ls i\ls t$, position $b_i$ is occupied and position $c_i$ is empty.
\item
$b_1=\beta$.
\item
For $1\ls i<t$, we have
\begin{itemize}
\item
$b_i\equiv c_i\ppmod e$, and all the positions $b_i-e,b_i-2e,\dots,c_i+e$ are occupied;
\item\label{empties}
all the positions $c_i-1,c_i-2,\dots,b_{i+1}+1$ are empty.
\end{itemize}
\item\label{mull1}
Either:
\begin{enumerate}
\item\label{mull1a}
$b_t\equiv c_t\ppmod e$, all the positions $b_t-e,\dots,c_t+e$ are occupied, and all the positions $c_t-1,c_t-2,\dots,\gamma$ are empty; or
\item\label{mull1b}
all the positions $b_t-e,b_t-2e,\dots$ are occupied and $c_t=\gamma$.
\end{enumerate}
\end{enumerate}

We define $\mu\dw$ to be the partition whose abacus display is obtained by moving the beads at positions $b_1,\dots,b_t$ to positions $c_1,\dots,c_t$, and we define the \emph{$e$-rim length} of $\mu$ to be $\erim(\mu) = |\mu|-|\mu\dw|=\sum_{i=1}^t(b_i-c_i)$.  It is straightforward to see that neither $\mu\dw$ nor $\erim(\mu)$ depends on the choice of $r$.
\end{defn}

\begin{eg}
Suppose $e=3$, and %take $r=$.
%Suppose 
$\mu=(12,11^2,7,6,5,3^2,2)$.  The abacus display for $\mu$ with $r=15$ is as follows.
\[\bsm
\tl&\tm&\tr\\
\bd&\bd&\bd\\
\bd&\bd&\bd\\
\nb&\nb&\bd\\
\nb&\bd&\bd\\
\nb&\nb&\bd\\
\nb&\bd&\nb\\
\bd&\nb&\nb\\
\nb&\nb&\bd\\
\bd&\nb&\bd\\
\nb&\nb&\nb
\esm\]	
We see that $\beta=26$ and $\gamma=6$.  We find that $t=3$ and
\[(b_1,c_1,b_2,c_2,b_3,c_3) = (26,20,18,15,14,6).\]
So $\erim(\mu) = 17$, and
\[\mu\dw = 
\bsm
\tl&\tm&\tr\\
\bd&\bd&\bd\\
\bd&\bd&\bd\\
\bd&\nb&\bd\\
\nb&\bd&\bd\\
\nb&\nb&\nb\\
\bd&\bd&\nb\\
\nb&\nb&\bd\\
\nb&\nb&\bd\\
\bd&\nb&\nb\\
\nb&\nb&\nb
\esm
=(10^2,8,5^2,2^2,1).\]
\end{eg}

Now we can describe the Mullineux map.  We define $m(\mu)$ recursively in $|\mu|$, setting $m(\varnothing)=\varnothing$.  If $\mu\neq\varnothing$, then we compute $\erim(\mu)$ and $\mu\dw$ as in Definition \ref{mulldef}.  Obviously $|\mu\dw|<|\mu|$, and we assume $m(\mu\dw)$ is defined.  Now set
\[l = \begin{cases}
\erim(\mu)-\mu'_1 & (e\mid \erim(\mu))\\
\erim(\mu)-\mu'_1+1 & (e\nmid \erim(\mu)),
\end{cases}\]
and define $m(\mu)$ to be the unique $e$-regular partition such that $(m(\mu))'_1=l$, $\erim(m(\mu))=\erim(\mu)$ and $(m(\mu))\dw = m(\mu\dw)$.  That this procedure always works (i.e.\ there is always a unique $e$-regular $m(\mu)$ with the required properties) is proved by Mullineux in \cite{mull}.

\subs{The Mullineux map and conjugation}

Given the statement of Theorem \ref{mainmull}, it will be helpful for us to study the map $\mu\mapsto m(\mu)'$ rather than $m$ itself.  To this end, we give a `conjugate' definition to Definition \ref{mulldef}.

\begin{defn}\label{mcdef}
Suppose $\nu$ is an $e$-restricted partition, and take an abacus display for $\nu$ with $r$ beads.  Let $\delta,\epsilon$ be the position of the last bead and the first empty space on the abacus, respectively.  Assuming $\nu\neq\varnothing$, there is a unique sequence $f_1>g_1>\dots>f_u>g_u$ of non-negative integers satisfying the following conditions.
\begin{enumerate}
\item
For each $1\ls i\ls u$, position $f_i$ is occupied and position $g_i$ is empty.
\item
$g_u=\epsilon$.
\item
For $1< i\ls u$, we have
\begin{itemize}
\item
$f_i\equiv g_i\ppmod e$, and all the positions $f_i-e,f_i-2e,\dots,g_i+e$ are empty;
\item
all the positions $f_i+1,f_i+2,\dots,g_{i-1}-1$ are occupied.
\end{itemize}
\item\label{mc1}
Either:
\begin{enumerate}
\item
$f_1\equiv g_1\ppmod e$, all the positions $f_1-e,\dots,g_1+e$ are empty, and all the positions $\delta,\delta-1,\dots,f_1+1$ are occupied; or
\item
all the positions $g_1+e,g_1+2e,\dots$ are empty and $f_1=\delta$.
\end{enumerate}
\end{enumerate}
We define $\nu\wod$ to be the partition whose abacus display is obtained by moving the beads at positions $f_1,\dots,f_u$ to positions $g_1,\dots,g_u$, and we define the \emph{conjugate $e$-rim length} of $\nu$ to be $\mire(\nu) = \sum_{i=1}^t(f_i-g_i)$.  It is straightforward to see that neither $\nu\wod$ nor $\mire(\nu)$ depends on the choice of $r$.
\end{defn}

Definition \ref{mcdef} is the result of applying Definition \ref{mulldef} to the $e$-regular partition $\nu'$, and then exploiting Lemma \ref{conbeta}.  This yields
\[\mire(\nu) = \erim(\nu')\]
and
\[\nu\wod = ((\nu')\dw)'.\]
Hence we can describe the map $\mu\mapsto m(\mu)'$, as follows.

\begin{lemma}\label{mc}
Suppose $\mu$ is an $e$-regular partition.  If $\mu=\varnothing$, then $m(\mu)'=\varnothing$.  Otherwise, set
\[l = \begin{cases}
\erim(\mu)-\mu'_1 & (e\mid \erim(\mu))\\
\erim(\mu)-\mu'_1+1 & (e\nmid \erim(\mu)).
\end{cases}\]
Then $m(\mu)'$ is the unique $e$-restricted partition such that $(m(\mu)')_1 = l$, $\mire(m(\mu)')=\erim(\mu)$ and $(m(\mu)')\wod = m(\mu\dw)'$.
\end{lemma}

We use Lemma \ref{mc} as our definition of the map $\mu\mapsto m(\mu)'$.

\begin{egs}
Suppose $e=3$, and take $r=6$.
\begin{enumerate}
\item
Suppose $\mu=(3)$ and $\nu=(2,1)$.  These partitions have the following abacus displays.
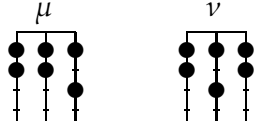
\[\begin{array}c
\mu\\
\bsm
\tl&\tm&\tr\\
\bd&\bd&\bd\\
\bd&\bd&\nb\\
\nb&\nb&\bd\\
\nb&\nb&\nb
\esm
\end{array}
\qquad
\begin{array}c
\nu\\
\bsm
\tl&\tm&\tr\\
\bd&\bd&\bd\\
\bd&\nb&\bd\\
\nb&\bd&\nb\\
\nb&\nb&\nb
\esm
\end{array}
\]
Applying Definition \ref{mulldef} to $\mu$, we see that $\beta=8$ and $\gamma=5$.  We have $t=1$, with $b_1=8$ and $c_1=5$.  So $\mu\dw=\varnothing$, $\erim(\mu)=3$ and $l=2$.

Applying Definition \ref{mcdef} to $\nu$, we have $u=1$, with $f_1=\delta=7$, $g_1=\epsilon=4$.  So $\nu\wod=\varnothing$, $\mire(\nu)=3$ and $\nu_1=2$, and we see that $m(\mu)'=\nu$.
\item
Now suppose $\mu=(6,3,1)$ and $\nu=(5,3,2)$.
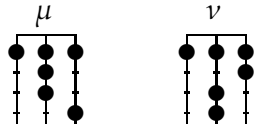
\[\begin{array}c
\mu\\
\bsm
\tl&\tm&\tr\\
\bd&\bd&\bd\\
\nb&\bd&\nb\\
\nb&\bd&\nb\\
\nb&\nb&\bd
\esm
\end{array}
\qquad
\begin{array}c
\nu\\
\bsm
\tl&\tm&\tr\\
\bd&\bd&\bd\\
\nb&\nb&\bd\\
\nb&\bd&\nb\\
\nb&\bd&\nb
\esm
\end{array}
\]
Applying Definition \ref{mulldef} to $\mu$, we have $\beta=11$, $\gamma=3$, $t=2$, and $(b_1,c_1,b_2,c_2)=(11,8,7,3)$.  So $\mu\dw=(3)$, $\erim(\mu)=7$, and $l=5$.

Applying Definition \ref{mcdef} to $\nu$, we have $\delta=10$, $\epsilon=3$, $u=1$, and $(f_1,g_1)=(10,3)$.  So $\nu\wod=(2,1)$ and $\mire(\nu)=7$.  By (1) above we have $m(\mu\dw)'=\nu\wod$, and since $\nu_1=5$ we have $m(\mu)'=\nu$.
\end{enumerate}
\end{egs}

\subs{Proof of Theorem \ref{mainmull}}

%\begin{notn}
%\subsubsection*{Notation}
Now we proceed with the proof of Theorem \ref{mainmull}; we begin with the following lemma.

\begin{lemma}\label{lltcor}
Suppose $\mu$ is an $e$-regular partition.  Then $\mu\dom m(\mu)'$.
\end{lemma}

\begin{pf}
This is immediate from Propositions \ref{lltthm} and \ref{qdecdom}.
\end{pf}

Now we fix some notation which will be in force for the remainder of Section \ref{mullsec}.  We fix an $e$-regular partition $\mu\neq\varnothing$ and set $\nu=m(\mu)'$.  We fix a large integer $r$, and let $\beta,\gamma, b_1,c_1,\dots,b_t,c_t$ be as in Definition \ref{mulldef}, and $\delta,\epsilon,f_1,g_1,\dots,f_u,g_u$ as in Definition \ref{mcdef}.  Let $y$ denote $\erim(\mu)$ ($=\mire(\nu)$).  Fix $k$ such that $\mu$ and $\nu$ are both $k$-empty, and let $c,d$ be as defined in \S\ref{absec}.  Let $x=d+(c-1)e$ be the position of the last bead on runner $d$ of the abacus display for $\mu$; since $\mu$ and $\nu$ have the same $e$-core (which is implicit in Lemma \ref{lltcor}), $x$ is also the position of the last bead on runner $d$ of the abacus display for $\nu$.
%\end{notn}

\begin{lemma}\label{lamustuff}
\indent
\vspace{-\topsep}
\begin{enumerate}
\item\label{lamustuff1}
We have
\[y = \begin{cases}
\delta-\gamma+1 & (e\mid y)\\
\delta-\gamma & (e\nmid y).\end{cases}\]
\item
None of $b_1,c_1,\dots,b_{t-1},c_{t-1}$ or $f_2,g_2,\dots,f_u,g_u$ is congruent to $d$ modulo $e$.
\item
$b_t-c_t$ is divisible by $e$ if and only if $f_1-g_1$ is.  If neither is divisible by $e$, then
\[
\gamma=c_t\equiv g_1\pmod{e}\qquad\text{and}
\qquad
\delta=f_1\equiv b_t\pmod{e}.\]
\end{enumerate}
\end{lemma}

\begin{pfenum}
\item
This follows from the statements $\gamma = \mu'_1-r$, $\delta=\nu_1+r-1$ and
\[\nu_1 = l = \begin{cases}
y-\mu'_1 & (e\mid y)\\
y-\mu'_1+1 & (e\nmid y).
\end{cases}\]
\item
Suppose $1\ls i<t$.  Then $b_i\equiv c_i\ppmod e$ and $b_i>c_i$.  But there is no bead on runner $d$ of the abacus display for $\mu$ with an empty space above it, so we cannot have $b_i\equiv d\ppmod e$.  Similarly $f_i,g_i\nequiv d$ if $1<i\ls u$.
\item
Since $y=\sum_i(b_i-c_i)$ and $b_i-c_i$ is divisible by $e$ for $i<t$, we have $b_t-c_t\equiv y\ppmod{e}$.  Similarly $f_1-g_1\equiv y\ppmod e$.  Now suppose neither $b_t-c_t$ nor $f_1-g_1$ is divisible by $e$.  Then in condition (\ref{mull1}) of Definition \ref{mulldef} we must be in case (b), and in condition (\ref{mc1}) of Definition \ref{mcdef} we must in be case (b).  So $c_t=\gamma$ and $f_1=\delta$.  Let $h,i,h',i'$ be the residues of $b_t,c_t,f_1,g_1$ modulo $e$.  Since $\mu$ and $\nu$ have the same $e$-core, there must be the same numbers of beads on corresponding runners of the abacus displays for $\mu$ and $\nu$.  Similarly, there are the same numbers of beads on corresponding runners of the abacus displays for $\mu\dw$ and $\nu\wod$.  To get from the abacus display of $\mu$ to the abacus display for $\mu\dw$, we move some beads up their runners, and then move a bead from runner $h$ to runner $i$.  Similarly, to get from $\nu$ to $\nu\wod$, we move some beads up their runners and then move a bead from runner $h'$ to runner $i'$.  Combining these statements, we see that $h=h'$ and $i=i'$.
\end{pfenum}

The proof of Theorem \ref{mainmull} is is by induction on $|\mu|$; the inductive step is to assume that the theorem holds with $\mu$ replaced by $\mu\dw$, and to compare $\ux k(\mu)$ with $\ux k(\mu\dw)$ and $\ux k(\nu)$ with $\ux k(\nu\wod)$.  The calculation required for this inductive step is broken into several parts.

\begin{lemma}\label{lamlem}
Suppose neither $f_1$ nor $g_1$ is congruent to $d$ modulo $e$.  Then $\nu\wod$ is $k$-empty, and
\[\ux k(\nu) = \ux k(\nu\wod)+\left\lfloor\frac{\delta-x}e\right\rfloor.\]
\end{lemma}

\begin{pf}
We obtain an abacus display for $\nu\wod$ by moving a bead from position $f_i$ to position $g_i$ for each $i$.  Since none of $f_1,g_1,\dots,f_u,g_u$ is congruent to $d$ modulo $e$, this has no effect on runner $d$ of the abacus display, and so $\nu\wod$ is $k$-empty.  Now we compute $\ux k(\nu) - \ux k(\nu\wod)$.

Suppose first that $\delta<x+e$; then we must show that $\ux k(\nu) - \ux k(\nu\wod)=0$.  For each $i$ we have $f_i\ls \delta<x+e$, so that (by Lemma \ref{altux}) moving a bead from position $f_i$ to position $g_i$ does not alter the value of $\ux k$, and we are done.

So we assume that $\delta>x+e$.  $g_u=\epsilon$ is the first empty position on the abacus display for $\nu$, and there is an empty space at position $x+e$, so we have $g_u\ls x+e$; but $g_u\nequiv d\ppmod e$ by assumption, so $g_u<x+e$.  Let $l$ be minimal such that $g_l<x+e$, and for $i=1,\dots,l$ write
\[g_i = d+a_ie+j_i\]
with $0<j_i<e$.  Also write
\[\delta = d+a_0e+j_0\]
with $0<j_0<e$.
\clam
For $1\ls i\ls l$ we have
\[d+a_{i-1}e < f_i < d+(a_{i-1}+1)e.\]
\prof
Suppose first that $i\gs2$.  Recall that every position between $f_i$ and $g_{i-1}$ in the abacus display for $\nu$ is occupied.  Since $g_{i-1}>x+e$, position $d+a_{i-1}e$ is unoccupied, so $d+a_{i-1}e$ does not lie between $f_i$ and $g_{i-1}$.  Hence
\[d+a_{i-1}e<f_i<g_{i-1}<d+(a_{i-1}+1)e,\]
as required.

Now consider $i=1$.  Either $f_1=\delta$ (in which case the result is immediate) or $f_1<\delta$ and every position between $f_1$ and $\delta$ in the abacus display for $\nu$ is occupied.  Assuming the latter and arguing as above, we get
\[d+a_0e<f_1<\delta<d+(a_0+1)e.\]
\malc

The claim, together with Lemma \ref{altux}, implies that if $1\ls i<l$ then moving a bead from position $f_i$ to position $g_i$ reduces the value of $\ux k$ by $a_{i-1}-a_i$.  Since $g_l<x+e$, moving the bead from position $f_l$ to position $g_l$ reduces the value of $\ux k$ by $a_{l-1}-(c-1)$.  For $l<i\ls t$, we have $f_i<g_l<x+e$, so moving a bead from position $f_i$ to position $g_i$ does not affect the value of $\ux k$.  Summing, we get
\begin{align*}
\ux k(\nu)-\ux k(\nu\wod) &= \sum_{i=1}^{l-1}(a_{i-1}-a_i)\ \ +\ \ (a_l-c+1)\\
&= a_0-c+1.
\end{align*}
On the other hand,
\begin{align*}
\left\lfloor\frac{\delta-x}e\right\rfloor &= \left\lfloor\frac{d+a_0 e+j_0 - d - (c-1)e}e\right\rfloor\\
&=a_0-c+1
\end{align*}
(since $0<j_0<e$), and we are done.
\end{pf}

\begin{lemma}\label{ctge}
Suppose $c_t>x+e$.  Then $y$ is divisible by $e$, $\mu\dw$ and $\nu\wod$ are both $k$-empty, and we have
\[\ux k(\mu)-\ux k(\mu\dw) = \frac ye,\qquad \ux k(\nu)-\ux k(\nu\wod)\gs \frac ye.\]
\end{lemma}

\begin{pf}
Since there is an empty space at position $x+e$ on the abacus display for $\mu$, we must have $\gamma\ls x+e<c_t$.  So in Definition \ref{mcdef}(\ref{mc1}) we must be in case (a), and hence $e\mid y$.  Since $b_i\equiv c_i$ and $f_i\equiv g_i\ppmod e$ for each $i$ and since it not possible to slide a bead up runner $d$ in the abacus display for either $\mu$ or $\nu$, we have $b_i,c_i,f_i,g_i\nequiv d\ppmod e$ for all $i$, so $\mu\dw$ and $\nu\wod$ are $k$-empty.

Now we examine $\mu$.  Since for each $i$ we have $x<c_i\equiv b_i\ppmod e$, Lemma \ref{altux} implies that moving the bead at position $b_i$ to position $c_i$ reduces the value of $\ux k$ by $(b_i-c_i)/e$.  So
\[\ux k(\mu) - \ux k(\mu\dw)=\frac{(b_1-c_1)+\dots+(b_t-c_t)}e =\frac ye\]
as required.

Next we examine $\nu$.  First we note that $\gamma>x$.  Indeed, Definition \ref{mulldef} tells us that every position between $c_t$ and $\gamma$ in the abacus for $\mu$ is empty; but position $x$ is occupied, so does not lie in this range.  Since $c_t>x$, we therefore have $\gamma>x$.

Now we have
\begin{align*}
\ux k(\nu)-\ux k(\nu\wod) &= \left\lfloor\frac{\delta-x}e\right\rfloor\tag*{by Lemma \ref{lamlem}}\\
&= \left\lfloor\frac{y+\gamma-x+1}e\right\rfloor\tag*{by Lemma \ref{lamustuff}(\ref{lamustuff1})}\\
&\gs \left\lfloor\frac ye\right\rfloor\\
&=\frac ye
\end{align*}
as required.
\end{pf}

\begin{lemma}\label{ctle}
Suppose $c_t<x+e$, and $b_t,c_t\nequiv d\ppmod e$.  Then $\mu\dw$ and $\nu\wod$ are $k$-empty, and we have
\[\ux k(\mu)-\ux k(\mu\dw) = \ux k(\nu)-\ux k(\nu\wod).\]
\end{lemma}

\begin{pfnb}
By Lemma \ref{lamlem}, we must show that
\[\ux k(\mu) - \ux k(\mu\dw) = \left\lfloor\frac{\delta-x}e\right\rfloor.\]
We use a calculation very similar to that used in the proof of Lemma \ref{lamlem}.  Let $l$ be maximal such that $b_l> x$; note that there is such an $l$, since $b_1=\beta\gs x$ and $b_1\nequiv x\ppmod e$.  For $l\ls i\ls t$, write
\[b_i = d+a_ie+j_i,\]
where $0<j_i<e$.
\clam
$c_l<x+e$.
\prof
If $l=t$ then this is true by assumption, so suppose $l<t$.  There is no bead on the abacus display for $\mu$ in any position between $c_l$ and $b_{l+1}$.  But there is a bead at position $x$, so $x$ does not lie between $c_l$ and $b_{l+1}$.  Since $b_{l+1}<x$, we have $c_l<x<x+e$.
\malc
The claim implies that moving a bead from position $b_l$ to position $c_l$ reduces the value of $\ux k$ by $a_l-(c-1)$.  For $1\ls i<l$, we have $c_i>b_l>x$, so moving a bead from position $b_i$ to position $c_i$ reduces the value of $\ux k$ by $(b_i-c_i)/e$.  For $l<i\ls t$, we have $b_i<x$, so moving a bead from position $b_i$ to position $c_i$ does not affect the value of $\ux k$.  So we have
\begin{align*}
\ux k(\mu)-\ux k(\mu\dw) &= \frac1e\sum_{i=1}^{l-1}(b_i-c_i)\quad+\quad a_l-(c-1)\\
&= \frac ye\quad-\quad\frac1e\sum_{i=l}^t(b_i-c_i) \quad+\quad a_l-(c-1).
\end{align*}
\clam
If $l\ls i<t$, then $c_i = d+a_{i+1}e+j_i$.
\prof
There are no beads in the abacus display for $\mu$ in any position between $c_i$ and $b_{i+1}$; since $b_{i+1}<x$, there is a bead at position $d+(a_{i+1}+1)e$, and therefore $d+(a_{i+1}+1)e$ does not lie between $c_i$ and $b_{i+1}$.  So we have
\[d+a_{i+1}e+j_{i+1}=b_{i+1}<c_i<d+(a_{i+1}+1)e;\]
since we know that $c_i\equiv b_i\ppmod e$, this implies that $j_i>j_{i+1}$ and $c_i=d+a_{i+1}e+j_i$.
\malc
Combining the claim with the expression above, we get
\begin{align*}
\ux k(\mu)-\ux k(\mu\dw) &= \frac ye\quad-\quad\sum_{i=l}^{t-1}(a_i-a_{i+1})\quad-\quad \frac1e(d+a_te+j_t-c_t)\quad +\quad a_l-(c-1)\\
&=\frac ye - (c-1) - \frac1e(d+j_t-c_t).
\end{align*}

Now we consider two cases, according to whether or not $e$ divides $y$.

\begin{description}
\item[\fbox{$e\nmid y$}]
Here we have $c_t=\gamma$ and $y=\delta-\gamma$, by Lemma \ref{lamustuff}.  So
\begin{align*}
\ux k(\mu)-\ux k(\mu\dw) &= \frac{\delta-\gamma-e(c-1)-d-j_t+\gamma}e\\
&= \frac{\delta-x-j_t}e\\
&=\left\lfloor\frac{\delta-x}e\right\rfloor,
\end{align*}
since $0<j_t<e$.
\item[\fbox{$e\mid y$}]
In this case, we write
\[\gamma = d+a^\ast e+j^\ast,\qquad c_t=d+a_\ast e+j_\ast\]
with $0<j^\ast,j_\ast<e$.  Since none of the positions $\gamma,\gamma+1,\dots,c_t$ is occupied but position $d+a_\ast e$ is occupied, we must have $a_\ast=a^\ast$,  so that $c_t=d+a^\ast e+j_t$.  Now the fact that $y=\delta-\gamma+1$ gives
\begin{align*}
\ux k(\mu)-\ux k(\mu\dw) &= \frac{\delta-\gamma+1-e(c-1)+a^\ast e}e\\
&= \frac{\delta-j^\ast+1-d-e(c-1)}e\\
&= \frac{\delta-x-(j^\ast-1)}e\\
&=\left\lfloor\frac{\delta-x}e\right\rfloor,
\end{align*}
since $0<j^\ast<e$.\hfill\qedsymbol
\end{description}\indent
\end{pfnb}

\begin{lemma}\label{dinv}
Suppose $b_t\equiv d\ppmod e$.  Then $\ux k(\mu)=0$.
\end{lemma}

\begin{pf}
First we note that, since it is impossible to move a bead up runner $d$, $c_t\nequiv d\ppmod e$, and this implies that $e\nmid y$.  So by Lemma \ref{lamustuff}(1) we have $y=\delta-\gamma$.  But $\delta=f_1$ is congruent to $d$ modulo $e$ (by Lemma \ref{lamustuff}(3)) and is the last occupied position on the abacus for $\nu$, and so must equal $x$.  So $y=x-\gamma$.

By Lemma \ref{altux}, the conclusion $\ux k(\mu)=0$ is the same as saying that $\beta<x+e$, so we prove the latter statement.  If $b_i<x$ for all $i$, then certainly $\beta=b_1<x+e$, so we assume otherwise, and let $l$ be maximal such that $b_l\gs x$.

For $l\ls i\ls t$ we write
\[b_i = d+a_ie+j_i,\]
with $0\ls j_i<e$.  Since $b_l,\dots,b_{t-1}\nequiv d\equiv b_t\ppmod e$, we actually have $j_t=0$ and $0<j_i<e$ for $l\ls i<t$.  Arguing as in the proof of Lemma \ref{ctle}, we have $j_i>j_{i+1}$ and $c_i = d+a_{i+1}e+j_i$ for $l\ls i<t$.  So
\begin{align*}
y\quad &=\quad \sum_{i=1}^t(b_i-c_i)\\
&=\quad \sum_{i=1}^{l-1}(b_i-c_i)\quad +\quad \sum_{i=l}^{t-1}(a_i-a_{i+1})e\quad +\quad d+a_te-\gamma\\
&=\quad \sum_{i=1}^{l-1}(b_i-c_i)\quad+\quad a_le+d-\gamma.
\end{align*}
Combining this with the equality $y=x-\gamma$ from above, we get
\[x\quad =\quad \sum_{i=1}^{l-1}(b_i-c_i)\quad +\quad a_le+d,\]
or
\[\sum_{i=1}^{l-1}(b_i-c_i)\quad+\quad b_l\quad =\quad x+j_l.\]
Now by assumption $b_l\gs x$, and so $\sum_{i=1}^{l-1}(b_i-c_i)\ls j_l<e$, which forces $l=1$.  And now we have
\[b_1 = x+j_1<x+e,\]
as required.
\end{pf}

\begin{pfof}{Theorem \ref{mainmull}}
Proceed by induction on $|\mu|$.  We consider several cases.
\begin{enumerate}
\item
First suppose $c_t>x+e$.  Then by Lemma \ref{ctge} $\mu\dw$ and $\nu\wod$ are $k$-empty, and
\[\ux k(\mu)-\ux k(\mu\dw) = \frac ye,\qquad \ux k(\nu)-\ux k(\nu\wod)\gs \frac ye.\]
By induction we have $\ux k(\mu\dw) = \ux k(\nu\wod)$, and so we get $\ux k(\mu)\ls \ux k(\nu)$.  By Lemmata \ref{lltcor} and \ref{uxdom} we have $\ux k(\mu)\gs\ux k(\nu)$, and the result follows.
\item
Next suppose that $c_t<x+e$, and that neither $b_t$ nor $c_t$ is congruent to $d$ modulo $e$.  Then by Lemma \ref{ctle} $\mu\dw$ and $\nu\wod$ are $k$-empty, and
\[\ux k(\mu)-\ux k(\mu\dw) = \ux k(\nu)-\ux k(\nu\wod),\]
and the result follows by induction.
\item
Next, suppose that $b_t\equiv f_1\equiv d\ppmod e$.  Then by Lemma \ref{dinv} we have $\ux k(\mu)=0$.  Since $\mu\dom\nu$ we have $\ux k(\mu)\gs \ux k(\nu)$, so $\ux k(\nu)=0$ too.
\item
Finally, consider the case where $c_t\equiv g_1\equiv d\ppmod e$.  Here, we replace the pair $(\mu,\nu)$ with $(\nu',\mu')$.  If we choose a large integer $s$ and let $\tilde b_1,\tilde c_1,\dots,\tilde b_{\tilde t},\tilde c_{\tilde t}$ be the integers given by Definition \ref{mulldef} with $\nu'$ in place of $\mu$ and $s$ in place of $r$, then by Lemma \ref{conbeta} we have $\tilde t=u$ and
\[\tilde b_i = r+s-1-g_{u+1-i},\qquad \tilde c_i = r+s-1-f_{u+1-i}\]
for each $i$.  If we set $\tilde d = (s+(e-1-k))\mmod e$, then we can compute $\tilde b_{\tilde t}\equiv \tilde d\ppmod e$; hence by case 3 above, we have $\ux{e-1-k}(\nu')=\ux{e-1-k}(\mu')$.  Now Corollary \ref{conux} gives the result.
\end{enumerate}
\end{pfof}

\end{document}